\documentclass[10.5pt]{article}
\usepackage[leqno]{amsmath}
\usepackage{amsfonts}
\usepackage{graphicx}

\usepackage{amsmath}
\usepackage{amssymb}
\usepackage{latexsym}
\usepackage{amsmath, amsfonts,amssymb, amsthm, euscript,makeidx,color,mathrsfs}

\def\sqr#1#2{{\vcenter{\vbox{\hrule height.#2pt
              \hbox{\vrule width.#2pt height#1pt \kern#1pt \vrule width.#2pt}
              \hrule height.#2pt}}}}
\def\signed #1{{\unskip\nobreak\hfil\penalty50
              \hskip2em\hbox{}\nobreak\hfil#1
              \parfillskip=0pt \finalhyphendemerits=0 \par}}
\def\endpf{\signed {$\sqr69$}}

\def\3n{\negthinspace \negthinspace \negthinspace }
\def\2n{\negthinspace \negthinspace }
\def\1n{\negthinspace }

\def\dbE{\mathbb{E}}
\def\dbF{\mathbb{F}}

\def\dbH{\mathbb{H}}

\def\dbP{\mathbb{P}}

\def\dbR{\mathbb{R}}
\def\dbS{\mathbb{S}}

\def\dbX{\mathbb{X}}
\def\dbY{\mathbb{Y}}
\def\dbZ{\mathbb{Z}}

\def\sQ{\mathscr{Q}}

\def\={\buildrel \triangle \over =}

\def\ds{\displaystyle}

\def\ns{\noalign{\ss}}
\def\e{\varepsilon}
\def\z{\zeta}

\def\l{\lambda}

\def\n{\nu}
\def\si{\sigma}

\def\f{\varphi}
\def\th{\theta}

\def\D{\Delta}
\def\Th{\Theta}
\def\L{\Lambda}

\def\O{\Omega}

\def\cF{{\cal F}}

\def\cR{{\cal R}}

\def\cU{{\cal U}}

\def\no{\noindent}

\def\ss{\smallskip}
\def\ms{\medskip}

\def\q{\quad}
\def\qq{\qquad}
\def\hb{\hbox}

\def\Ra{\mathop{\Rightarrow}}

\def\lan{\mathop{\langle}}
\def\ran{\mathop{\rangle}}

\def\esssup{\mathop{\rm esssup}}

\def\h{\widehat}
\def\wt{\widetilde}

\def\cd{\cdot}

\def\ae{\hbox{\rm a.e.{ }}}
\def\as{\hbox{\rm a.s.{ }}}

\def\tr{\hbox{\rm tr$\,$}}

\def\deq{\mathop{\buildrel\D\over=}}
\def\les{\leqslant}
\def\ges{\geqslant}

\def\({\Big (}
\def\){\Big )}
\def\[{\Big[}
\def\]{\Big]}
\def\bde{\begin{definition}}
\def\ede{\end{definition}}
\def\be{\begin{equation}}
\def\bel{\begin{equation}\label}
\def\ee{\end{equation}}
\def\bt{\begin{theorem}}
\def\et{\end{theorem}}
\def\bc{\begin{corollary}}
\def\ec{\end{corollary}}
\def\bl{\begin{lemma}}
\def\el{\end{lemma}}
\def\bp{\begin{proposition}}
\def\ep{\end{proposition}}
\def\bas{\begin{assumption}}
\def\eas{\end{assumption}}
\def\br{\begin{remark}}
\def\er{\end{remark}}
\def\ba{\begin{array}}
\def\ea{\end{array}}
\def\ed{\end{document}}

\def\square#1{\vbox{\hrule\hbox{\vrule height#1%
     \kern#1\vrule}\hrule}}
\def\rectangle#1#2{\vbox{\hrule\hbox{\vrule height#1%
     \kern#2\vrule}\hrule}}

\font\tenbb=msbm10 \font\sevenbb=msbm7 \font\fivebb=msbm5

\newfam\bbfam
\scriptscriptfont\bbfam=\fivebb \textfont\bbfam=\tenbb
\scriptfont\bbfam=\sevenbb

\newtheorem{lemma}{Lemma}[section]
\newtheorem{remark}{Remark}[section]

\newtheorem{theorem}{Theorem}[section]
\newtheorem{corollary}{Corollary}[section]

\newtheorem{definition}{Definition}[section]
\newtheorem{proposition}{Proposition}[section]
\newtheorem{assumption}{Assumption}[section]

\makeatletter
   
   \@addtoreset{equation}{section}
\makeatother

\begin{document}

\title{\bf Linear Quadratic Stochastic Differential Games: \\
Open-Loop and Closed-Loop Saddle Points\thanks{This work is
supported in part by NSF Grant DMS-1007514, and China Scholarship
Council.}}
\author{Jingrui Sun\thanks{School of Mathematical Sciences, University of
Science and Technology of China, Hefei, Anhui, 230026, P. R. China;
Email: sjr@mail.ustc.edu.cn.}\ , \ \ and \ \ Jiongmin Yong
\thanks{Department of Mathematics, University of Central Florida, Orlando, FL 32816, USA;
Email: jiongmin.yong@ucf.edu.}
 }
\maketitle

\noindent {\bf Abstract:} In this paper, we consider a linear
quadratic stochastic two-person zero-sum
differential game. The controls for both players are allowed to
appear in both drift and diffusion of the state equation. The
weighting matrices in the performance functional are not assumed to
be definite/non-singular. A necessary and sufficient condition for
the existence of a closed-loop saddle point is established in terms
of the solvability of a Riccati differential equation with {\it
certain regularity}. It is possible that the closed-loop saddle
point fails to exist, and at the same time, the corresponding
Riccati equation admits a solution (which does not have needed
regularity). Also, we will indicate that the solution of the Riccati
equation may be non-unique. \\


\ms

\no\bf AMS Mathematics Subject Classification. \rm 93E20, 91A23,
49N70.

\section{Introduction}

Let $(\O,\cF,\dbF,\dbP)$ be a given complete filtered probability
space along with a one-dimensional standard Brownian motion
$W=\{W(t),\cF_t; 0\les t < \infty \}$, where $\dbF=\{\cF_t\}_{t\ges0}$
is the natural filtration of $W$ augmented by all the $\dbP$-null
sets in $\cF$ \cite{Karatzas-Shreve 1991,Yong-Zhou 1999}. Consider
the following controlled linear stochastic differential equation on
$[t,T]$:
\bel{state}\left\{\2n\ba{ll}
\ns\ds dX(s)=\big[A(s)X(s)+B_1(s)u_1(s)+B_2(s)u_2(s)+b(s)\big]ds\\
\ns\ds\qq\qq+\big[C(s)X(s)+D_1(s)u_1(s)+D_2(s)u_2(s)+\si(s)\big]dW(s),\q s\in[t,T], \\
\ns\ds X(t)= x.\ea\right.\ee
In the above, $X(\cd)$ is called the {\it state process} taking
values in the $n$-dimensional Euclidean space $\dbR^n$ with the {\it
initial state} $x$ at the initial time $t$; for $i=1,2$, $u_i(\cd)$
is called the {\it control process} of Player $i$ taking values in
$\dbR^{m_i}$. We assume that $A(\cd), B_1(\cd), B_2(\cd), C(\cd),
D_1(\cd)$, $D_2(\cd)$ are deterministic matrix-valued functions of proper
dimensions, and $b(\cd)$, $\si(\cd)$ are vector-valued $\dbF$-adapted processes.
For any $t\in[0,T)$, we define
$$\cU_i[t,T]=\Big\{u_i:[t,T]\times\O\to\dbR^{m_i}\bigm|u_i(\cd)\hb{ is $\dbF$-adapted, }
\dbE\int_t^T|u_i(s)|^2ds<\infty\Big\},\qq i=1,2.$$
Any element $u_i(\cd)\in\cU_i[t,T]$ is called an {\it admissible
control} of Player $i$ on $[t,T]$. Under some mild conditions on the
coefficients, for any {\it initial pair} $(t,x)\in[0,T)\times\dbR^n$
and control pair $(u_1(\cd),u_2(\cd))\in\cU_1[t,T]\times\cU_2[t,T]$,
state equation (\ref{state}) admits a unique solution $X(\cd)\equiv
X(\cd\,;t,x,u_1(\cd),u_2(\cd))$. To measure the performance of the
controls $u_1(\cd)$ and $u_2(\cd)$, we introduce the following
functional:
\bel{cost}\ba{ll}
\ns\ds J(t,x;u_1(\cd),u_2(\cd))\deq{1\over2}\dbE\Big\{\lan
GX(T),X(T)\ran+2\lan g,X(T)\ran\\
\ns\ds\qq+\int_t^T\Big[\lan\begin{pmatrix}Q(s)&S_1(s)^T&S_2(s)^T\\
                                          S_1(s)&R_{11}(s)&R_{12}(s)\\
                                          S_2(s)&R_{21}(s)&R_{22}(s)\end{pmatrix}
                                          \begin{pmatrix}X(s)\\ u_1(s)\\ u_2(s)
                                          \end{pmatrix},\begin{pmatrix}X(s)\\ u_1(s)\\ u_2(s)
                                          \end{pmatrix}\ran
+2\lan\begin{pmatrix}q(s)\\ \rho_1(s)\\ \rho_2(s)\end{pmatrix},
                                          \begin{pmatrix}X(s)\\ u_1(s)\\ u_2(s)
                                          \end{pmatrix}\ran\Big]ds\Big\},\ea\ee
where $Q(\cd)$, $S_1(\cd)$, $S_2(\cd)$, $R_{11}(\cd)$,
$R_{12}(\cd)$, $R_{21}(\cd)$, $R_{22}(\cd)$ are deterministic
matrix-valued functions of proper dimensions with $Q(\cd)^T=Q(\cd)$,
$R_{ii}(\cd)^T=R_{ii}(\cd)$ ($i=1,2$), $R_{12}(\cd)^T=R_{21}(\cd)$,
and $G$ is a symmetric matrix; $q(\cd)$, $\rho_1(\cd)$,
$\rho_2(\cd)$ are allowed to be vector-valued $\dbF$-adapted
processes, and $g$ is allowed to be an $\cF_T$-measurable random
variable. We assume that (\ref{cost}) is a cost functional for
Player 1, and a payoff functional for Player 2. Therefore, Player 1
wishes to minimize (\ref{cost}) by selecting a control process
$u_1(\cd)\in\cU_1[t,T]$, while Player 2 wishes to maximize
(\ref{cost}) by selecting a control process $u_2(\cd)\in\cU_2[t,T]$.
The above described problem is referred to as a linear quadratic
(LQ, for short) stochastic {\it two-person zero-sum differential
game}, denoted by Problem (SG). The study of LQ deterministic
two-person zero-sum differential games, denoted by Problem (DG),
which is a special case of Problem (SG) where the stochastic part is
absent, can be traced back to the work of Ho--Bryson--Baron
\cite{Ho-Bryson-Baron 1965} in 1965. In 1970, Schmitendorf studied
both open-loop and closed-loop strategies for Problem (DG)
(\cite{Schmitendorf 1970}). Among other things, it was shown that
the existence of a closed-loop saddle point may not imply that of an
open-loop saddle point. In 1979, Bernhard carefully investigated
Problem (DG) from closed-loop point of view (\cite{Bernhard 1979});
see also the book by Basar and Bernhard \cite{Basar-Bernhard 1991}
in this aspect. In 2005, Zhang \cite{Zhang 2004} proved that for a
Problem (DG), the existence of the open-loop value is equivalent to
the finiteness of the corresponding open-loop lower and upper
values, which is also equivalent to the existence of an open-loop
saddle point. Along this line, there were a couple of follow-up
works \cite{Delfour 2007,Delfour-Sbarba 2009} appeared afterwards.
In 2006, Mou--Yong studied a Problem (SG) from an open-loop point of
view by means of Hilbert space method (\cite{Mou-Yong 2006}). The
main purpose of this paper is to study Problem (SG) from both
open-loop and closed-loop point of view.

\ms

We know that when $m_2=0$, Problem (SG) is reduced to a linear
quadratic stochastic optimal control problem, denoted by Problem
(SLQ). Thus, Problem (SLQ) can be regarded as a special case of
Problem (SG). See \cite{Chen-Yong 2000,Ait Rami-Moore-Zhou 2001,
Chen-Yong 2001,Hu-Zhou 2003,Tang 2003,McAsey-Mou 2006,Qian-Zhou
2013} for some relevant results on Problem (SLQ). Further, linear
quadratic deterministic optimal control problem, denoted by Problem
(DLQ), can be regarded as a special case of Problem (SLQ) and
Problem (DG). The history of Problem (DLQ) can further be traced
back to the work of Bellman--Glicksberg--Gross
(\cite{Belman-Gicksberg-Gross 1958}) in 1958, Kalman (\cite{Kalman
1960}), and Letov (\cite{Letov 1960}) in 1960. See \cite{Yong-Zhou
1999} for some historic remarks on Problems (DLQ) and (SLQ).

\ms

The rest of the paper is organized as follows. Section 2 will
collect some preliminary results. Among other things, we will state
some results for Problem (SLQ), which is a special case of our main
results presented in later sections. In Section 3, we pose our
differential game problem, carefully explaining the open-loop and
closed-loop saddle points. Section 4 is devoted to the study of
open-loop saddle points by variational method. In Section 5, we
characterize closed-loop saddle points by means of Riccati equation.
Some examples are presented in Section 6, and several concluding
remarks are collected in Section 7.

\section{Preliminaries}

We recall that $\dbR^n$ is the $n$-dimensional Euclidean space,
$\dbR^{n\times m}$ is the space of all $(n\times m)$ matrices,
endowed with the inner product $(M,N)\mapsto\tr[M^TN]$, and
$\dbS^n\subseteq \dbR^{n\times n}$ is the set of all $(n\times n)$
symmetric matrices. We recall that for any $M\in\dbR^{m\times n}$,
there exists a unique matrix $M^\dag\in\dbR^{n\times m}$, called the
(Moore-Penrose) {\it pseudo-inverse} of $M$, satisfying the
following (\cite{Penrose 1955}):
$$MM^\dag M=M,\q M^\dag MM^\dag=M^\dag,\q(MM^\dag)^T=MM^\dag,\q
(M^\dag M)^T=M^\dag M.$$
In addition, if $M=M^T\in\dbS^n$, then
$$\ba{ll}
\ns\ds M^\dag=(M^\dag)^T,\q MM^\dag=M^\dag M;\qq\hb{and}\qq
M\ges0\iff M^\dag\ges0.\ea$$

Next, let $T>0$ be a fixed time horizon. For any $t\in[0,T)$ and
Euclidean space $\dbH$, let
$$\ba{ll}
\ns\ds C([t,T];\dbH)=\Big\{\f:[t,T]\to\dbH\bigm|\f(\cd)\hb{ is
continuous }\Big\},\\
\ns\ds
L^p(t,T;\dbH)=\Big\{\f:[t,T]\to\dbH\bigm|\int_t^T|\f(s)|^pds<\infty\Big\},\q1\les p<\infty,\\
\ns\ds
L^\infty(t,T;\dbH)=\Big\{\f:[t,T]\to\dbH\bigm|\esssup_{s\in[t,T]}|\f(s)|<\infty\Big\}.\ea$$
Next, we recall that $(\O,\cF,\dbF,\dbP)$ is a complete filtered
probability space and $W(\cd)$ is a one-dimensional standard
Brownian motion with $\dbF=\{\cF_t\}_{t\ges0}$ being its natural
filtration augmented by all the $\dbP$-null sets in $\cF$. We denote
$$\ba{ll}
\ns\ds L^2_{\cF_T}(\O;\dbH)=\Big\{\xi:\O\to\dbH\bigm|\f(\cd)\hb{ is
$\cF_T$-measurable, }\dbE|\xi|^2<\infty\Big\},\\
\ns\ds
L_\dbF^2(t,T;\dbH)=\Big\{\f:[t,T]\times\O\to\dbH\bigm|\f(\cd)\hb{ is
$\dbF$-adapted, }\dbE\int^T_t|\f(s)|^2ds<\infty\Big\},\\
\ns\ds
L_\dbF^2(\Omega;C([t,T];\dbH))=\Big\{\f:[t,T]\times\O\to\dbH\bigm|\f(\cd)\hb{
is $\dbF$-adapted, continuous, }\dbE\big(\sup_{t\les s\les
T}|\f(s)|^2\big)<\infty\Big\},\\
\ns\ds L^2_\dbF(\O;L^1(0,T;\dbR^n))=\Big\{\f:[0,T]\times
\O\to\dbR^n\ \bigm|\ \f(\cd)\ \hb{is $\dbF$-adapted,
}\dbE\Big(\int_0^T|\f(t)|dt\Big)^2<\infty\Big\}.\ea$$

\ms

In the rest of this section, we look at Problem (SLQ). We write the
corresponding state equation as follows:
\bel{SDE}\left\{\2n\ba{ll}
\ns\ds
dX(s)=\big[A(s)X(s)+B(s)u(s)+b(s)\big]ds+\big[C(s)X(s)+D(s)u(s)+\si(s)\big]dW(s),\\
\ns\ds\qq\qq\qq\qq\qq\qq\qq\qq\qq\qq\qq\qq\qq s\in[t,T],\\
\ns\ds X(t)=x.\ea\right.\ee
The cost functional takes the following form:
\bel{cost2}\ba{ll}
\ns\ds J(t,x;u(\cd))\deq{1\over2}\dbE\Big\{\lan GX(T),X(T)\ran+2\lan
g,X(T)\ran\\
\ns\ds\qq\qq\qq\q+\int_t^T\Big[\lan\begin{pmatrix}Q(s)&S(s)^T\\
                                                  S(s)&R(s)\end{pmatrix}
                                   \begin{pmatrix}X(s)\\
                                                  u(s)\end{pmatrix},
                                   \begin{pmatrix}X(s)\\
                                                  u(s)\end{pmatrix}\ran
                                                  +2\lan\begin{pmatrix}q(s)\\
                                                    \rho(s)\end{pmatrix},
                                   \begin{pmatrix}X(s)\\
                                                  u(s)\end{pmatrix}\ran\Big]ds\Big\}.\ea\ee
We adopt the following assumptions.

\ms

{\bf(S1)} The coefficients of the state equation satisfy the
following:
$$\left\{\2n\ba{ll}
\ns\ds A(\cd)\in L^1(0,T;\dbR^{n\times n}),\q B(\cd)\in
L^2(0,T;\dbR^{n\times m}),\q b(\cd)\in L^2_\dbF(\O;L^1(0,T;\dbR^n)),\\
\ns\ds C(\cd)\in L^2(0,T;\dbR^{n\times n}),\q D(\cd)\in
L^\infty(0,T;\dbR^{n\times m}),\q\si(\cd)\in
L_\dbF^2(0,T;\dbR^n).\ea\right.$$

{\bf(S2)} The weighting coefficients in the cost functional satisfy
the following:
$$\left\{\2n\ba{ll}
\ns\ds Q(\cd)\in L^1(0,T;\dbS^n),\q S(\cd)\in L^2(0,T;\dbR^{m\times
n}),\q R(\cd)\in L^\infty(0,T;\dbS^m),\\
\ns\ds q(\cd)\in L^2_\dbF(\O;L^1(0,T;\dbR^n)),\q\rho(\cd)\in
L_\dbF^2(0,T;\dbR^m),\q G\in\dbS^n,\q g\in
L^2_{\cF_T}(\O;\dbR^n).\ea\right.$$

We note that under (S1), for any $(t,x)\in[0,T)\times\dbR^n$ and
$u(\cd)\in\cU[t,T]\equiv L^2_\dbF(t,T;\dbR^m)$, the state equation
(\ref{SDE}) admits a unique solution $X(\cd)\equiv
X(\cd\,;t,x,u(\cd))$. Further, under (S2), the cost functional
(\ref{cost2}) is well-defined for every $(t,x)\in[0,T)\times\dbR^n$
and $u(\cd)\in\cU[t,T]$. Therefore, the following problem is
meaningful.

\ms

\bf Problem (SLQ). \rm For any given initial pair
$(t,x)\in[0,T)\times\dbR^n$, find a $\bar u(\cd)\in\cU[t,T]$ such
that
\bel{2.43}V(t,x)\deq J(t,x;\bar
u(\cd))=\inf_{u(\cd)\in\cU[t,T]}J(t,x;u(\cd))\les
J(t,x;u(\cd)),\qq\forall u(\cd)\in\cU[t,T].\ee

\ms

Any $\bar u(\cd)\in\cU[t,T]$ satisfying (\ref{2.43}) is called an
{\it open-loop optimal control} of Problem (SLQ) for $(t,x)$, the
corresponding $\bar X(\cd)\equiv X(\cd\,;t,x,\bar u(\cd))$ is called
an {\it open-loop optimal state process} and $(\bar X(\cd),\bar
u(\cd))$ is called an {\it open-loop optimal pair}. The function
$V(\cd\,,\cd)$ is called the {\it value function} of Problem (SLQ).

\ms

The following result is concerned with open-loop optimal control of
Problem (SLQ) for a given initial pair.

\ms

\bf Theorem 2.1. \sl Let {\rm(S1)--(S2)} hold. For a given initial
pair $(t,x)\in[0,T)\times\dbR^n$, a state-control pair $(\bar
X(\cd),\bar u(\cd))$ is an open-loop optimal pair of Problem
{\rm(SLQ)} if and only if the following {\it stationarity condition}
holds:
\bel{2.8}B(s)^T\bar Y(s)+D(s)^T\bar Z(s)+S(s)\bar X(s)+R(s)\bar
u(s)+\rho(s)=0,\qq \ae s\in[t,T],~\as\ee
where $(\bar Y(\cd),\bar Z(\cd))$ is the {\it adapted solution} to
the following {\it backward stochastic differential equation} (BSDE,
for short):
\bel{BSDE0}\left\{\2n\ba{ll}
\ns\ds d\bar Y(s)=-\big[A(s)^T\bar Y(s)+C(s)^T\bar Z(s)+Q(s)\bar
X(s)+S(s)^T\bar
u(s)+q(s)\big]+\bar Z(s)dW(s),\\
\ns\ds\qq\qq\qq\qq\qq\qq\qq\qq\qq\qq\qq\qq\qq\qq\qq s\in[t,T],\\
\ns\ds\bar Y(T)=G\bar X(T)+g,\ea\right.\ee
and the following {\it convexity
condition} holds: For any $u(\cd)\in\cU[t,T]$,
\bel{>00}\dbE\Big\{\1n\lan
GX_0(T),X_0(T)\ran\1n+\2n\int_t^T\1n\2n\big[\lan
Q(s)X_0(s),X_0(s)\ran\1n+2\lan S(s)X_0(s),u(s)\ran\1n+\1n\lan
R(s)u(s),u(s)\ran\big]ds\Big\}\ges0,\ee
where $X_0(\cd)$ is the solution to the following:
\bel{X0s}\left\{\2n\ba{ll}
\ns\ds dX_0(s)=\big[A(s)X_0(s)+B(s)u(s)\big]ds+\big[C(s)X_0(s)+D(s)u(s)\big]dW(s),\qq s\in[t,T],\\
\ns\ds X_0(t)=0.\ea\right.\ee

\it Proof. \rm Suppose $(\bar X(\cd),\bar u(\cd))$ is  a
state-control pair corresponding to the given initial pair
$(t,x)\in[0,T)\times\dbR^n$. For any $u(\cd)\in\cU[t,T]$ and $\e\in\dbR$,
let $X^\e(\cd)=X(\cd\,;t,x,\bar u(\cd)+\e u(\cd))$. Then
$$\left\{\2n\ba{ll}
\ns\ds dX^\e(s)=\big\{A(s)X^\e(s)+B(s)\big[\bar u(s)+\e
u(s)\big]+b(s)\big\}ds\\
\ns\ds\qq\qq\q+\big\{C(s)X^\e(s)+D(s)\big[\bar u(s)+\e u(s)\big]+\si(s)\big\}dW(s),\qq s\in[t,T],\\
\ns\ds X^\e(t)=x.\ea\right.$$
Thus, $X_0(\cd)\equiv{X^\e(\cd)-\bar X(\cd)\over\e}$ is independent
of $\e$ and satisfies (\ref{X0s}). Then
$$\ba{ll}
\ns\ds J(t,x;\bar u(\cd)+\e u(\cd))-J(t,x;\bar u(\cd))\\
\ns\ds={\e\over2}\dbE\Big\{\lan G\big[2\bar X(T)+\e
X_0(T)\big],X_0(T)\ran
+2\lan g,X_0(T)\ran\\
\ns\ds\qq+\int_t^T\Big[\lan\begin{pmatrix}Q(s)&S(s)^T\\
                                                  S(s)&R(s)\end{pmatrix}
                                   \begin{pmatrix}2\bar X(s)+\e X_0(s)\\
                                                  2\bar u(s)+\e u(s)\end{pmatrix},
                                   \begin{pmatrix}X_0(s)\\
                                                  u(s)\end{pmatrix}\ran
                                                  +2\lan\begin{pmatrix}q(s)\\
                                                    \rho(s)\end{pmatrix},
                                   \begin{pmatrix}X_0(s)\\
                                                  u(s)\end{pmatrix}\ran\Big]ds\Big\}\\
\ns\ds=\e\dbE\Big\{\lan G\bar X(T),X_0(T)\ran+\lan g,X_0(T)\ran\\
\ns\ds\qq+\int_t^T\Big[\lan Q\bar X,X_0\ran+\lan S\bar X,u\ran+\lan
SX_0,\bar u\ran+\lan R\bar u,u\ran+\lan
q,X_0\ran+\lan\rho,u\ran\Big]ds\Big\}\\
\ns\ds\qq+{\e^2\over2}\dbE\Big\{\lan
GX_0(T),X_0(T)\ran+\int_t^T\Big[\lan
Q(s)X_0(s),X_0(s)\ran+2\lan S(s)X_0(s),u(s)\ran\\
\ns\ds\qq\qq\qq\qq\qq\qq\qq\qq\qq+\lan
R(s)u(s),u(s)\ran\Big]ds\Big\}.\ea$$
Now, let $(\bar Y(\cd),\bar Z(\cd))$ be the adapted solution to the
BSDE (\ref{BSDE0}). Then
$$\ba{ll}
\ns\ds\dbE\Big\{\lan G\bar X(T)+g,X_0(T)\ran+\int_t^T\Big[\lan
Q\bar X+S^T\bar u+q,X_0\ran+\lan S\bar X+R\bar
u+\rho,u\ran\Big]ds\Big\}\\
\ns\ds\q=\dbE\Big\{\int_t^T\Big[\lan-(A^T\bar Y+C^T\bar Z+Q\bar
X+S^T\bar u+q),X_0\ran+\lan\bar Y,AX_0+Bu\ran+\lan
\bar Z,CX_0+Du\ran\\
\ns\ds\qq\qq\qq+\lan Q\bar X+S^T\bar u+q,X_0\ran+\lan S\bar X+R\bar
u+\rho,u\ran\Big]ds\Big\}\\
\ns\ds\q=\dbE\int_t^T\lan B^T\bar Y+D^T\bar Z+S\bar X+R\bar
u+\rho,u\ran ds.\ea$$
Hence,
$$\ba{ll}
\ns\ds J(t,x;\bar u(\cd)+\e u(\cd))=J(t,x;\bar
u(\cd))+\e\dbE\Big\{\int_t^T\lan B^T\bar Y+D^T\bar Z+S\bar X+R\bar
u+\rho,u\ran ds\Big\}\\
\ns\ds\qq+{\e^2\over2}\dbE\Big\{\lan
GX_0(T),X_0(T)\ran+\int_t^T\Big[\lan
Q(s)X_0(s),X_0(s)\ran+2\lan S(s)X_0(s),u(s)\ran\\
\ns\ds\qq\qq\qq\qq\qq\qq\qq\qq\qq+\lan
R(s)u(s),u(s)\ran\Big]ds\Big\}.\ea$$
Therefore, $(\bar X(\cd),\bar u(\cd))$ is an open-loop optimal
control of Problem (SLQ) for $(t,x)$ if and only if (\ref{2.8}) and
(\ref{>00}) hold. \endpf

\ms

Next, for any $t\in[0,T)$, take $\Th(\cd)\in L^2(t,T;\dbR^{m\times
n})\equiv\sQ[t,T]$, and $v(\cd)\in\cU[t,T]$. For any $x\in\dbR^n$,
we consider the following equation:
\bel{closed-loop1}\left\{\2n\ba{ll}
\ns\ds
dX(s)=\big\{\big[A(s)+B(s)\Th(s)\big]X(s)+B(s)v(s)+b(s)\big\}ds\\
\ns\ds\qq\qq+\big\{\big[C(s)+D(s)\Th(s)\big]X(s)+D(s)v(s)+\si(s)\big\}dW(s),\q
s\in[t,T],\\
\ns\ds X(t)=x,\ea\right.\ee
which admits a unique solution $X(\cd)\equiv
X(\cd\,;t,x,\Th(\cd),v(\cd))$, depending on $\Th(\cd)$ and $v(\cd)$.
The above is called a {\it closed-loop system} of the original state
equation (\ref{SDE}) under {\it closed-loop strategy}
$(\Th(\cd),v(\cd))$. We point out that $(\Th(\cd),v(\cd))$ is
independent of the initial state $x$. With the above solution
$X(\cd)$, we define
$$\ba{ll}
\ns\ds J(t,x;\Th(\cd)X(\cd)+v(\cd))={1\over2}\dbE\Big\{\lan
GX(T),X(T)\ran+2\lan g,X(T)\ran\\
\ns\ds\qq\qq\qq\qq\qq\qq+\int_t^T\Big[\lan\begin{pmatrix}Q(s)&S(s)^T\\
                                                  S(s)&R(s)\end{pmatrix}
                                   \begin{pmatrix}X(s)\\
                                                  \Th(s)X(s)+v(s)\end{pmatrix},
                                   \begin{pmatrix}X(s)\\
                                                  \Th(s)+v(s)\end{pmatrix}\ran\\
\ns\ds\qq\qq\qq\qq\qq\qq\qq+2\lan\begin{pmatrix}q(s)\\
                                                    \rho(s)\end{pmatrix},
                                   \begin{pmatrix}X(s)\\
                                                  \Th(s)+v(s)\end{pmatrix}\ran\Big]ds\Big\}.\ea$$

We now introduce the following definition.

\ms

\bf Definition 2.2. \rm A pair $(\bar\Th(\cd),\bar
v(\cd))\in\sQ[t,T]\times\cU[t,T]$ is called a {\it closed-loop
optimal strategy} of Problem (SLQ) on $[t,T]$ if
\bel{2.45}J(t,x;\bar\Th(\cd)\bar X(\cd)+\bar v(\cd))\les
J(t,x;u(\cd)),\qq\forall x\in\dbR^n,\q u(\cd)\in\cU[t,T],\ee
where $\bar X(\cd)=X(\cd\,;t,x,\bar\Th(\cd),\bar v(\cd))$.

\ms

\bf Remark 2.3. \rm In the above, both $\bar\Th(\cd)$ and $\bar
v(\cd)$ are required to be independent of the initial state
$x\in\dbR^n$. Also, for fixed initial pair
$(t,x)\in[0,T)\times\dbR^n$, (\ref{2.45}) implies that the outcome
$$\bar u(\cd)\equiv\bar\Th(\cd)\bar X(\cd)+\bar v(\cd)\in\cU[t,T]$$
of the closed-loop optimal strategy $(\bar\Th(\cd),\bar v(\cd))$ is
an open-loop optimal control of Problem (SLQ) for $(t,x)$.
Therefore, for Problem (SLQ), the existence of closed-loop optimal
strategy on $[t,T]$ implies the existence of open-loop optimal
control for initial pair $(t,x)$ for any $x\in\dbR^n$.

\ms

For closed-loop optimal strategies, we have the following
characterization.

\ms

\bf Theorem 2.4. \sl Let {\rm(S1)--(S2)} hold. Then Problem
{\rm(SLQ)} admits a closed-loop optimal strategy if and only if the
following Riccati equation admits a solution $P(\cd)\in C([t,T];\dbS^n)$:
\bel{Riccati2.52}\left\{\2n\ba{ll}
\ns\ds\dot P+PA+A^TP+C^TPC+Q\\
\ns\ds\qq-(PB+C^TPD+S^T)(R+D^TPD)^\dag(B^TP+D^TPC+S)=0,\qq\ae s\in[t,T],\\
\ns\ds\cR\big(B^TP+D^TPC+S\big)\subseteq\cR\big(R+D^TPD\big),\qq\ae s\in[t,T],\\
\ns\ds R+D^TPD\ges0,\qq\ae s\in[t,T],\\
\ns\ds P(T)=G,\ea\right.\ee
such that
\bel{2.53}\big[R(\cd)+D(\cd)^TP(\cd)D(\cd)\big]^\dag
\big[B(\cd)^TP(\cd)+D(\cd)^TP(\cd)C(\cd)+S(\cd)\big]\in
L^2(t,T;\dbR^{m\times n}),\ee
and the adapted solution $(\eta(\cd),\z(\cd))$ of the following BSDE:
\bel{BODE1}\left\{\2n\ba{ll}
\ns\ds d\eta=-\Big\{\big[A^T-(PB+C^TPD+S^T)(R+D^TPD)^\dag
B^T\big]\eta\\
\ns\ds\qq\qq+\big[C^T-(PB+C^TPD+S^T)(R+D^TPD)^\dag D^T\big]\z\\
\ns\ds\qq\qq+\big[C^T-(PB+C^TPD+S^T)(R+D^TPD)^\dag
D^T\big]P\si\\
\ns\ds\qq\qq-(PB+C^TPD+S^T)(R+D^TPD)^\dag
\rho+Pb+q\Big\}ds+\z dW(s),\\
\ns\ds B^T\eta+D^T\z+D^TP\si+\rho\in\cR(R+D^TPD),\q \ae s\in[t,T],~\as\\
\ns\ds\eta(T)=g,\ea\right.\ee
satisfies
\bel{}\big[R(\cd)+D(\cd)^TP(\cd)D(\cd)\big]^\dag\big[B(\cd)^T\eta(\cd)+D(\cd)^T\z(\cd)
+D(\cd)^TP(\cd)\si(\cd)+\rho(\cd)\big]\in L_\dbF^2(t,T;\dbR^m).\ee
In this case, any closed-loop optimal strategy $(\bar\Th(\cd),\bar
v(\cd))$ of Problem {\rm(SLQ)} admits the following representation:
\bel{Th-v2.54}\left\{\2n\ba{ll}
\ns\ds\bar\Th\1n=\1n-(R+\1n D^T\1n PD)^\dag(B^T\1n P\1n+\1n D^T\1n PC\1n+\1n S)\1n
+\1n\big[I\1n-(R\1n+\1n D^T\1n PD)^\dag(R\1n+\1n D^T\1n PD)\big]\th,\\
\ns\ds\bar v=-(R+D^T\1n PD)^\dag(B^T\1n\eta+\1nD^T\1n\z+\1n D^T\1n
P\si\1n+\1n\rho)+\big[I-(R\1n+\1n D^T\1n PD)^\dag(R+\1n D^T\1n PD)
\big]\n,\ea\right.\ee
for some $\th(\cd)\in L^2(t,T;\dbR^{m\times n})$ and $\n(\cd)\in
L_\dbF^2(t,T;\dbR^m)$. Further, the value function admits the
following representation:
\bel{Value2.18}\ba{ll}
\ns\ds V(t,x)\equiv\inf_{u(\cd)\in\cU[t,T]}J(t,x;u(\cd))\\
\ns\ds\qq\q={1\over2}\,\dbE\Big\{\lan
P(t)x,x\ran+2\lan\eta(t),x\ran+\1n\int_t^T\3n\big[\lan
P\si,\si\ran\1n+\1n2\lan\eta,b\ran+\1n2\lan\z,\si\ran\\
\ns\ds\qq\qq\qq\q-\1n \lan(R\1n+\1n D^T\1n
PD)^\dag(B^T\1n\eta\1n+D^T\1n\z\1n+\1n D^T\1n P\si\1n+\1n\rho),B^T\1n\eta\1n+D^T\1n\z\1n+\1n
D^T\1n P\si\1n+\1n\rho\ran\big]ds\Big\}.\ea\ee

\ms

\rm

We prefer not to give a proof of the above result here since it will
be a special case of a similar result for Problem (SG), which will be
presented below.

\section{Stochastic Differential Games}

We now return to our Problem (SG). Recall the sets
$\cU_i[t,T]=L^2_\dbF(t,T;\dbR^{m_i})$ of all open-loop controls of
Player $i$ ($i=1,2$). For notational simplicity, we let $m=m_1+m_2$
and denote
$$\ba{ll}
\ns\ds B(\cd)=(B_1(\cd),B_2(\cd)),\q D(\cd)=(D_1(\cd),D_2(\cd)),\\
\ns\ds S(\cd)=\begin{pmatrix}S_1(\cd)\\ S_2(\cd)\end{pmatrix},\q R(\cd)=\begin{pmatrix}R_{11}(\cd)&R_{12}(\cd)\\
R_{21}(\cd)&R_{22}(\cd)\end{pmatrix}\equiv\begin{pmatrix}R_1(\cd)\\ R_2(\cd)\end{pmatrix},\q\rho(\cd)=\begin{pmatrix}\rho_1(\cd)\\
\rho_2(\cd)\end{pmatrix},\q u(\cd)=\begin{pmatrix}u_1(\cd)\\
u_2(\cd)\end{pmatrix}.\ea$$
Naturally, we identify $\cU[t,T]=\cU_1[t,T]\times\cU_2[t,T]$. With
such notations, the state equation becomes
\bel{state3.1}\left\{\2n\ba{ll}
\ns\ds
dX(s)=\big[A(s)X(s)+B(s)u(s)+b(s)\big]ds+\big[C(s)X(s)+D(s)u(s)+\si(s)\big]dW(s),
\q s\in[t,T],\\
\ns\ds X(t)=x,\ea\right.\ee
and the performance functional becomes
\bel{cost3.2}\ba{ll}
\ns\ds
J(t,x;u_1(\cd),u_2(\cd))=J(t,x;u(\cd))={1\over2}\,\dbE\Big\{\lan
GX(T),X(T)\ran+2\lan g,X(T)\ran\\
\ns\ds\qq\qq\qq+\int_t^T\Big[\lan\begin{pmatrix}Q(s)&S(s)^T\\
                                          S(s)&R(s)\end{pmatrix}\begin{pmatrix}
                                          X(s)\\ u(s)\end{pmatrix},
                                          \begin{pmatrix}
                                          X(s)\\
                                          u(s)\end{pmatrix}\ran+
                                          2\lan\begin{pmatrix}q(s)\\ \rho(s)
                                          \end{pmatrix},\begin{pmatrix}X(s)\\
                                          u(s)\
                                          \end{pmatrix}\ran\Big] ds\Big\}.\ea\ee
When $b(\cd),\si(\cd),q(\cd),\rho(\cd),g(\cd)=0$, we denote the
problem by Problem $\hb{(SG)}^0$, which is a special case of Problem
(SG). With the above notation, we introduce the following standard
assumptions:

\ms

{\bf(SG1)} The coefficients of the state equation satisfy the
following:
$$\left\{\2n\ba{ll}
\ns\ds A(\cd)\in L^1(0,T;\dbR^{n\times n}),\q B(\cd)\in
L^2(0,T;\dbR^{n\times m}),\q b(\cd)\in L^2_\dbF(\O;L^1(0,T;\dbR^n)),\\
\ns\ds C(\cd)\in L^2(0,T;\dbR^{n\times n}),\q D(\cd)\in
L^\infty(0,T;\dbR^{n\times m}),\q\si(\cd)\in
L_\dbF^2(0,T;\dbR^n).\ea\right.$$

{\bf(SG2)} The weighting coefficients in the cost functional satisfy
the following:
$$\left\{\2n\ba{ll}
\ns\ds Q(\cd)\in L^1(0,T;\dbS^n),\q S(\cd)\in L^2(0,T;\dbR^{m\times
n}),\q R(\cd)\in L^\infty(0,T;\dbS^m),\\
\ns\ds q(\cd)\in L^2_\dbF(\O;L^1(0,T;\dbR^n)),\q\rho(\cd)\in
L_\dbF^2(0,T;\dbR^m),\q G\in\dbS^n,\q g\in
L^2_{\cF_T}(\O;\dbR^n).\ea\right.$$

Under (SG1), for any $(t,x)\in[0,T)\times\dbR^n$, and
$u(\cd)\equiv(u_1(\cd)^T,u_2(\cd)^T)^T\in\cU_1[t,T]\times\cU_2[t,T]
\equiv\cU[t,T]$, equation (\ref{state3.1}) admits a unique solution
(\cite{Yong-Zhou 1999})
$$X(\cd)\deq X(\cd\,;t,x,u_1(\cd),u_2(\cd))\equiv X(\cd\,;t,x,u(\cd))\in L^2_\dbF\big(\O;C([0,T];\dbR^n)\big).$$
Moreover, the following estimate holds:
$$\dbE\(\sup_{t \les s\les T}|X(s)|^2\)\les
K\dbE\Big\{|x|^2+\(\int_t^T|b(s)|ds\)^2+\int_t^T|\si(s)|^2ds
+\int^T_t|u(s)|^2ds\Big\},$$
hereafter, $K>0$ represents a generic constant which can be
different from line to line. Therefore, under (SG1)--(SG2), the
quadratic performance functional $J(t,x;u(\cd))\equiv
J(t,x;u_1(\cd),u_2(\cd))$ is well defined for all
$(t,x)\in[0,T)\times\dbR^n$ and
$(u_1(\cd),u_2(\cd))\in\cU_1[t,T]\times\cU_2[t,T]$. Having the
above, we now introduce the following definition.

\ms

\bf Definition 3.1. \rm (i) A pair
$(u^*_1(\cd),u^*_2(\cd))\in\cU_1[t,T]\times\cU_2[t,T]$ is called an
{\it open-loop saddle point} of Problem (SG) for the initial pair
$(t,x)\in[0,T)\times\dbR^n$ if for any
$(u_1(\cd),u_2(\cd))\in\cU_1[t,T]\times\cU_2[t,T]$,
\bel{saddle-open}J(t,x;u^*_1(\cd),u_2(\cd))\les
J(t,x;u^*_1(\cd),u^*_2(\cd))\les J(x;u_1(\cd),u^*_2(\cd)).\ee

(ii) The {\it open-loop upper value} $V^+(t,x)$ of Problem (SG) at
$(t,x)\in[0,T)\times\dbR^n$ and the {\it open-loop lower value}
$V^-(t,x)$ of Problem (SG) at $(t,x)\in[0,T)\times\dbR^n$ are
defined by the following:
\bel{}\left\{\ba{ll}
\ns\ds
V^+(t,x)=\inf_{u_1(\cd)\in\cU_1[t,T]}\sup_{u_2(\cd)\in\cU_2[t,T]}
J(t,x;u_1(\cd),u_2(\cd)),\\
\ns\ds
V^-(t,x)=\sup_{u_2(\cd)\in\cU_2[t,T]}\inf_{u_1(\cd)\in\cU_1[t,T]}
J(t,x;u_1(\cd),u_2(\cd)),\ea\right.\ee
which automatically satisfy the following:
$$V^-(t,x)\les V^+(t,x),\qq(t,x)\in[0,T)\times\dbR^n.$$
In the case that
\bel{}V^-(t,x)=V^+(t,x)\equiv V(t,x),\ee
we say that Problem (SG) admits an {\it open-loop value} $V(t,x)$ at
$(t,x)$. The maps $(t,x)\mapsto V^\pm(t,x)$ and $(t,x)\mapsto
V(t,x)$ are called {\it open-loop upper value function}, {\it
open-loop lower value function}, and {\it open-loop value function},
respectively.

\ms

Next, we let
$$\sQ_i[t,T]=L^2(t,T;\dbR^{m_i\times n}),\qq i=1,2.$$
For any initial pair $(t,x)\in[0,T)\times\dbR^n$,
$\Th(\cd)\equiv(\Th_1(\cd)^T,\Th_2(\cd)^T)^T\in\sQ_1[t,T]\times\sQ_2[t,T]$
and
$v(\cd)\equiv(v_1(\cd)^T,v_2(\cd)^T)^T\in\cU_1[t,T]\times\cU_2[t,T]$,
consider the following system:
\bel{state-closed}\left\{\2n\ba{ll}
\ns\ds dX(s)=\Big\{\big[A(s)+B(s)\Th(s)\big]X(s)+B(s)v(s)+b(s)\Big\}ds\\
\ns\ds\qq\qq\q+\Big\{\big[C(s)+D(s)\Th(s)\big]X(s)+D(s)v(s)+\si(s)\Big\}dW(s),
\qq s\in[t,T], \\
\ns\ds X(t)= x.\ea\right.\ee
Clearly, under (SG1), the above admits a unique solution
$X(\cd)\equiv X(\cd\,;t,x,\Th_1(\cd),v_1(\cd);\Th_2(\cd),v_2(\cd))$.
If we denote
$$u_i(\cd)=\Th_i(\cd)X(\cd)+v_i(\cd),\qq i=1,2,$$
then the above (\ref{state-closed}) coincides with the original
state equation (\ref{state}). We refer to (\ref{state-closed}) as
a {\it closed-loop system} of the original system. With the solution
$X(\cd)$ to (\ref{state-closed}), we denote
$$\ba{ll}
\ns\ds J\big(t,x;\Th_1(\cd)X(\cd)+v_1(\cd),\Th_2(\cd)X(\cd)+v_2(\cd)\big)
\equiv J(t,x;\Th(\cd)X(\cd)+v(\cd))\\
\ns\ds={1\over2}\,\dbE\Big\{\lan
GX(T),X(T)\ran+2\lan g,X(T)\ran\\
\ns\ds\q+\int_t^T\Big[\lan\begin{pmatrix}Q(s)&S(s)^T\\
                                          S(s)&R(s)\end{pmatrix}\begin{pmatrix}
                                          X(s)\\ \Th(s)X(s)+v(s)\end{pmatrix},
                                          \begin{pmatrix}
                                           X(s)\\ \Th(s)X(s)+v(s)\end{pmatrix}
                                           \ran\\
\ns\ds\qq\qq+2\lan\begin{pmatrix}q(s)\\ \rho(s)
                  \end{pmatrix},\begin{pmatrix}X(s)\\ \Th(s)X(s)+v(s)\end{pmatrix}
                  \ran\Big]ds\Big\}\\
\ns\ds={1\over2}\,\dbE\Big\{\lan
GX(T),X(T)\ran+2\lan g,X(T)\ran\\
\ns\ds\qq+\1n\int_t^T\2n\Big[\lan\begin{pmatrix}Q\1n+\1n\Th^T\1n
S\1n+\1n S^T\1n\Th\1n+\1n\Th^T\1n R\Th&
                                          S^T\2n+\1n\Th^T\1n R\\
                                          S+R\Th&R\end{pmatrix}\begin{pmatrix}
                                          X\\ v\end{pmatrix},\begin{pmatrix}X\\ v
                                          \end{pmatrix}\ran\1n
                                          +2\lan\begin{pmatrix}q\1n+\1n\Th^T\1n\rho\\
\rho\end{pmatrix},\begin{pmatrix}X\\
v\end{pmatrix}\ran\Big]ds\Big\}.\ea$$
Similarly, one can define
$J(t,x;\Th_1(\cd)X(\cd)+v_1(\cd),u_2(\cd))$,
$J(t,x;u_1(\cd),\Th_2(\cd)X(\cd)+v_2(\cd))$. We now introduce the
following definition.

\ms

\bf Definition 3.2. \rm A 4-tuple
$(\Th_1^*(\cd),v_1^*(\cd);\Th_2^*(\cd),v_2^*(\cd))\in
\sQ_1[t,T]\times\cU_1[t,T]\times\sQ_2[t,T]\times\cU_2[t,T]$ is
called a {\it closed-loop saddle point} of Problem (SG) on $[t,T]$
if for any $x\in\dbR^n$ and
$(u_1(\cd),u_2(\cd))\in\cU_1[t,T]\times\cU_2[t,T]$,
\bel{saddle-closed}\ba{ll}
\ns\ds J(t,x;\Th_1^*(\cd)X(\cd)+v_1^*(\cd),u_2(\cd))\les
J(t,x;\Th_1^*(\cd)X^*(\cd)+v_1^*(\cd),\Th_2^*(\cd)X^*(\cd)+v_2^*(\cd))\\
\ns\ds\qq\qq\qq\qq\qq\qq\q~\les
J(t,x;u_1(\cd),\Th_2^*(\cd)X(\cd)+v_2^*(\cd)).\ea\ee

\ms

There are some remarks in order.

\ms

(i) An open-loop saddle point $(u_1^*(\cd),u_2^*(\cd))$ usually
depends on the initial state $x$, whereas, a closed-loop saddle
point $(\Th_1^*(\cd),v_1^*(\cd);\Th_2^*(\cd),v_2^*(\cd))$ is
required to be independent of the initial state $x$.

\ms

(ii) In (\ref{saddle-closed}), the state process $X(\cd)$ appearing
in $J(t,x;\Th^*_1(\cd)X(\cd)+v_1^*(\cd),u_2(\cd))$ is different from
that in $J(t,x;u_1(\cd),\Th^*_2(\cd)X(\cd)+v_2^*(\cd))$; and both
are different from $X^*(\cd)\equiv
X(\cd\,;t,x,\Th_1^*(\cd),v_1^*(\cd);\Th_2^*(\cd),v_2^*(\cd))$ which
is the solution of (\ref{state-closed}) corresponding to
$$(\Th_1(\cd),v_1(\cd);\Th_2(\cd),v_2(\cd))=(\Th_1^*(\cd),v_1^*(\cd);
\Th_2^*(\cd),v_2^*(\cd)).$$
Therefore, comparing with (\ref{saddle-open}), we see that
(\ref{saddle-closed}) does not imply that
$(\Th_1^*(\cd)X^*(\cd)+v_1^*(\cd),\Th_2^*(\cd)X^*(\cd)+v_2^*(\cd))$
is an open-loop saddle point of Problem (SG), for the initial pair
$(t,X^*(t))$. This is different from Problem (SLQ) (for which the
outcome $\bar\Th(\cd)\bar X(\cd)+\bar v(\cd)$ of a closed-loop
optimal strategy $(\bar\Th(\cd),\bar v(\cd))$ is an open-loop
optimal control for the initial pair $(t,\bar X(t))$).

\ms

More precisely, let us comparing the following two inequalities:
\bel{open3.10}J(t,x;u_1^*(\cd),u_2^*(\cd))\les
J(t,x;u_1(\cd),u_2^*(\cd)),\ee
and
\bel{closed3.11}J(t,x;\Th_1^*(\cd)X^*(\cd)+v_1^*(\cd),\Th_2^*(\cd)X^*(\cd)+v_2^*(\cd))\les
J(t,x;u_1(\cd),\Th_2^*(\cd)X(\cd)+v_2^*(\cd)).\ee
For (\ref{open3.10}), we look at the following state equation:
$$\left\{\2n\ba{ll}
\ns\ds dX(s)=\big[A(s)X(s)+B_1(s)u_1(s)+B_2(s)u_2^*(s)+b(s)\big]ds\\
\ns\ds\qq\qq\q+\big[C(s)X(s)+D_1(s)u_1(s)+D_2(s)u_2^*(s)+\si(s)\big]dW(s),\q
s\in[t,T],\\
\ns\ds X(t)=x,\ea\right.$$
and the following cost functional
$$\ba{ll}
\ns\ds J_1(t,x;u_1(\cd))\equiv
J(t,x;u_1(\cd),u_2^*(\cd))={1\over2}\dbE\Big\{\lan
GX(T),X(T)\ran+2\lan g,X(T)\ran\\
\ns\ds\qq\qq+\int_t^T\Big[\lan QX,X\ran+2\lan S_1X,u_1\ran+\lan
R_{11}u_1,u_1\ran+\lan R_{22}u_2^*,u_2^*\ran+2\lan R_{12}u_2^*,u_1\ran\\
\ns\ds\qq\qq\qq\qq+2\lan S_2X,u_2^*\ran+2\lan
q,X\ran+2\lan\rho_1,u_1\ran+2\lan\rho_2,u_2^*\ran\Big]ds\Big\}\\
\ns\ds={1\over2}\dbE\Big\{\lan GX(T),X(T)\ran+2\lan
g,X(T)\ran+\int_t^T\Big[\lan QX,X\ran+2\lan S_1X,u_1\ran+\lan
R_{11}u_1,u_1\ran\\
\ns\ds\qq\qq+2\lan
q+S_2^Tu_2^*,X\ran+2\lan\rho_1+R_{12}u_2^*,u_1\ran+\lan
R_{22}u_2^*,u_2^*\ran+2\lan\rho_2,u_2^*\ran\Big]ds\Big\}.\ea$$
Therefore, (\ref{open3.10}) holds if and only if $u_1^*(\cd)$ is an
open-loop optimal control of Problem (SLQ) with
\bel{3.12}\left\{\2n\ba{ll}
\ns\ds\wt A=A,\q\wt B=B_1,\q\wt b=b+B_2u_2^*,\q\wt C=C,\q\wt D=D_1,\q\wt\si=\si+D_2u_2^*,\\
\ns\ds\wt G=G,\q\wt g=g,\q\wt Q=Q,\q\wt S=S_1,\q\wt R=R_{11},\q\wt
q=q+S_2^Tu_2^*,\q\wt\rho=\rho_1+R_{12}u_2^*.\ea\right.\ee
However, for (\ref{closed3.11}), we look at the following state
equation:
$$\left\{\2n\ba{ll}
\ns\ds
dX_1(s)=\Big\{\big[A(s)+B_2(s)\Th_2^*(s)\big]X_1(s)+B_1(s)u_1(s)+B_2(s)v_2^*(s)+b(s)
\Big\}ds\\
\ns\ds\qq\qq\q+\Big\{\big[C(s)\1n+\1n
D_2(s)\Th_2^*(s)\big]X_1(s)\1n+\1n D_1(s)u_1(s)\1n+\1n
D_2(s)v_2^*(s)\1n
+\1n\si(s)\Big\}dW(s),\q s\in[t,T],\\
\ns\ds X_1(t)=x,\ea\right.$$
and the following cost functional
$$\ba{ll}
\ns\ds
\bar J_1(t,x;u_1(\cd))=J(t,x;u_1(\cd),\Th_2^*(\cd)X_1(\cd)+v_2^*(\cd))\\
\ns\ds={1\over2}\dbE\Big\{\lan GX_1(T),X_1(T)\ran+2\lan
g,X_1(T)\ran+\int_t^T\Big[\lan QX_1,X_1\ran+\lan
R_{11}u_1,u_1\ran\\
\ns\ds\qq\qq+\lan
R_{22}(\Th_2^*X_1+v_2^*),\Th_2^*X_1+v_2^*\ran+2\lan
S_1X_1,u_1\ran+2\lan S_2X_1,\Th_2^*X_1+v_2^*\ran\\
\ns\ds\qq\qq+2\lan R_{21}u_1,\Th_2^*X_1+v_2^*\ran+2\lan
q,X_1\ran+2\lan\rho_1,u_1\ran
+2\lan\rho_2,\Th_2^*X_1+v_2^*\ran\Big]ds\Big\}\\
\ns\ds={1\over2}\dbE\Big\{\lan GX_1(T),X_1(T)\ran+2\lan
g,X_1(T)\ran+\int_t^T\Big[\lan QX_1,X_1\ran+\lan
R_{11}u_1,u_1\ran\\
\ns\ds\qq\qq+\lan(\Th_2^*)^TR_{22}\Th_2^*X_1,X_1\ran+2\lan
(\Th_2^*)^TR_{22}v_2^*,X_1\ran+\lan R_{22}v_2^*,v_2^*\ran+2\lan
S_1X_1,u_1\ran\\
\ns\ds\qq\qq+\lan[S_2^T\Th_2^*+(\Th_2^*)^TS_2]X_1,X_1\ran +2\lan
S_2^Tv_2^*,X_1\ran+2\lan R_{12}\Th_2^*X_1,u_1\ran+2\lan
R_{12}v_2^*,u_1\ran\\
\ns\ds\qq\qq+2\lan
q,X_1\ran+2\lan\rho_1,u_1\ran+2\lan(\Th_2^*)^T\rho_2,X_1\ran+2\lan\rho_2,
v_2^*\ran\Big]ds\Big\}\\
\ns\ds={1\over2}\dbE\Big\{\lan GX_1(T),X_1(T)\ran+2\lan
g,X_1(T)\ran\\
\ns\ds\qq\qq+\int_t^T\Big[\lan[Q+(\Th_2^*)^TR_{22}\Th_2^*+(\Th_2^*)^TS_2
+S_2^T\Th_2^*]X_1,X_1\ran+\lan
R_{11}u_1,u_1\ran\\
\ns\ds\qq\qq+2\lan(S_1+R_{12}\Th_2^*)X_1,u_1\ran
+2\lan q+[S_2^T+(\Th_2^*)^TR_{22}]v_2^*+(\Th_2^*)^T\rho_2,X_1\ran\\
\ns\ds\qq\qq+2\lan\rho_1+R_{12}v_2^*,u_1\ran+\lan
R_{22}v_2^*,v_2^*\ran+2\lan\rho_2,v_2^*\ran\Big]ds\Big\}.\ea$$
Then, $(\Th_1^*(\cd),v_1^*(\cd))$ is a closed-loop optimal strategy
for a Problem (SLQ), with
\bel{3.13}\left\{\2n\ba{ll}
\ns\ds\wt A=A+B_2\Th_2^*,\q\wt B=B_1,\q\wt b=b+B_2v_2^*,\\
\ns\ds\wt C=C+D_2\Th_2^*,\q\wt D=D_1,\q\wt\si=\si+D_2v_2^*,\\
\ns\ds\wt
Q=Q+(\Th_2^*)^TR_{22}\Th_2^*+(\Th_2^*)^TS_2+S_2^T\Th_2^*,\q\wt
S=S_1+R_{12}\Th_2^*,\q\wt R=R_{11},\\
\ns\ds\wt q=q+[S_2^T+(\Th_2^*)^TR_{22}]v_2^*+(\Th_2^*)^T\rho_2,\q
\wt\rho=\rho_1+R_{12}v_2^*,\\
\ns\ds\wt G=G,\q\wt g=g.\ea\right.\ee
Comparing (\ref{3.12}) and (\ref{3.13}), we see that one cannot say
anything whether the outcome $\Th_1^*(\cd)X^*(\cd)+v_1^*(\cd)$ of
$(\Th_1^*(\cd),v_1^*(\cd))$ for the initial pair $(t,x)$ has
anything to do with $u_1^*(\cd)$.

\ms

On the other hand, the following result, which is similar to
Berkovitz's equivalence lemma for Problem (DG) found in
\cite{Berkovitz 1971}, will tells something a little differently and
will be useful below.

\ms

\bf Proposition 3.3. \sl Let {\rm(SG1)--(SG2)} hold. For
$(\Th_i^*(\cd),v_i^*(\cd))\in\sQ_i[t,T]\times\cU_i[t,T]$, the
following statements are equivalent:

\ms

{\rm(i)} $(\Th_1^*(\cd),v_1^*(\cd);\Th_2^*(\cd),v_2^*(\cd))$ is a
closed-loop saddle point of Problem {\rm(SG)} on $[t,T]$.

\ms

{\rm(ii)} For any $x\in\dbR^n$,
$(\Th_1(\cd),\Th_2(\cd))\in\sQ_1[t,T]\times\sQ_2[t,T]$ and
$(v_1(\cd),v_2(\cd))\in\cU_1[t,T]\times\cU_2[t,T]$, the following
holds:
\bel{}\ba{ll}
\ns\ds
J(t,x;\Th^*_1(\cd)X(\cd)+v_1^*(\cd),\Th_2(\cd)X(\cd)+v_2(\cd))\les
J(t,x;\Th^*_1(\cd)X^*(\cd)+v_1^*(\cd),\Th^*_2(\cd)X^*(\cd)+v_2^*(\cd))\\
\ns\ds\q\qq\qq\qq\qq\qq\qq\qq\qq\qq\les
J(t,x;\Th_1(\cd)X(\cd)+v_1(\cd),\Th^*_2(\cd)X(\cd)+v_2^*(\cd)).\ea\ee

{\rm(iii)} For any $x\in\dbR^n$ and $(v_1(\cd),v_2(\cd))\in
\cU_1[t,T]\times\cU_2[t,T]$, the following holds:
\bel{}\ba{ll}
\ns\ds
J(t,x;\Th^*_1(\cd)X(\cd)+v_1^*(\cd),\Th^*_2(\cd)X(\cd)+v_2(\cd))\les
J(t,x;\Th^*_1(\cd)X^*(\cd)+v_1^*(\cd),\Th^*_2(\cd)X^*(\cd)+v_2^*(\cd))\\
\ns\ds\qq\qq\qq\qq\qq\qq\qq\qq\qq\q\les
J(t,x;\Th^*_1(\cd)X(\cd)+v_1(\cd),\Th^*_2(\cd)X(\cd)+v_2^*(\cd)).\ea\ee

\ms

\it Proof. \rm (i) $\Ra$ (ii). For any $\Th_i(\cd)\in\sQ_i[t,T]$ and
$v_i(\cd)\in\cU_i[t,T]$, $i=1,2$, let $X(\cd)$ be the solution to
the following SDE:
\bel{}\left\{\2n\ba{ll}
\ns\ds
dX(s)=\Big\{\big[A(s)+B_1(s)\Th_1(s)+B_2(s)\Th_2^*(s)\big]X(s)+B_1(s)v_1(s)
+B_2(s)v^*_2(s)+b(s)\Big\}ds \\
\ns\ds\qq\qq+\Big\{\big[C(s)\1n+\1n D_1(s)\Th_1(s)\1n+\1n
D_2(s)\Th^*_2(s)\big]X(s)\1n+\1n D_1(s)v_1(s)\1n+\1n
D_2(s)v^*_2(s)\1n+\1n\si(s)
\Big\}dW(s),\\
\ns\ds\qq\qq\qq\qq\qq\qq\qq\qq\qq\qq\qq\qq\qq\qq\qq s\1n\in\1n[t,T],\\
\ns\ds X(t)=x.\ea\right.\ee
Set
$$u_1(\cd)\deq\Th_1(\cd)X(\cd)+v_1(\cd)\in\cU_1[t,T].$$
By uniqueness, $X(\cd)$ also solves the following SDE:
\bel{}\left\{\2n\ba{ll}
\ns\ds
dX(s)=\Big\{\big[A(s)+B_2(s)\Th^*_2(s)\big]X(s)+B_1(s)u_1(s)+B_2(s)v^*_2(s)+b(s)
\Big\}dt\\
\ns\ds\qq\qq+\Big\{\big[C(s)\1n+\1n D_2(s)\Th^*_2(s)\big]X(s)\1n+\1n
D_1(s)u_1(s)\1n+\1n D_2(s)v^*_2(s)\1n+\1n\si(s)\Big\}dW(s),
\q s\in[t,T],\\
\ns\ds X(t)=x.\ea\right.\ee
Therefore,
$$\ba{ll}
\ns\ds
J(t,x;\Th^*_1(\cd)X^*(\cd)+v_1^*(\cd),\Th^*_2(\cd)X^*(\cd)+v^*_2(\cd))
\les J(t,x;u_1(\cd),\Th^*_2(\cd)X(\cd)+v_2^*(\cd))\\
\ns\ds=J(t,x;\Th_1(\cd)X(\cd)+v_1(\cd),\Th^*_2(\cd)X(\cd)+v_2^*(\cd)).\ea$$
Similarly, we have
$$J(t,x;\Th^*_1(\cd)X(\cd)+v_1^*(\cd),\Th_2(\cd)X(\cd)+v_2(\cd))\les
J(t,x;\Th^*_1(\cd)X^*(\cd)+v^*_1(\cd),\Th^*_2(\cd)X^*(\cd)+v_2^*(\cd)).$$

\ms

(ii) $\Ra$ (iii) is trivial, by taking $\Th_i(\cd)=\Th^*_i(\cd)$,
$i=1,2$.

\ms

(iii) $\Ra$ (i). For any $x\in\dbR^n$, and any
$u_1(\cd)\in\cU_1[t,T]$, let $X(\cd)$ be the solution of the
following SDE:
\bel{}\left\{\2n\ba{ll}
\ns\ds dX(s)=\Big\{\big[A(s)+B_2(s)\Th^*_2(s)\big]X(s)+B_1(s)u_1(s)+B_2(s)v^*_2(s)+b(s)\Big\}ds\\
\ns\ds\qq\qq\q+\Big\{\big[C(s)+D_2(s)\Th^*_2(s)\big]X(s)+D_1(s)u_1(s)
+D_2(s)v^*_2(s)+\si(s)\Big\}dW(s),\\
\ns\ds\qq\qq\qq\qq\qq\qq\qq\qq\qq\qq\qq\qq\qq\qq s\in[t,T],\\
\ns\ds X(t)=x.\ea\right.\ee
Set
$$v_1(\cd)=u_1(\cd)-\Th^*_1(\cd)X(\cd)\in\cU_1[t,T],$$
then $X(\cd)$ is also the solution to the following SDE:
\bel{}\left\{\2n\ba{ll}
\ns\ds
dX(s)=\Big\{\big[A(s)+B_1(s)\Th^*_1(s)+B_2(s)\Th^*_2(s)\big]X(s)
+B_1(s)v_1(s)+B_2(s)v^*_2(s)+b(s)\Big\}ds \\
\ns\ds\qq\qq+\Big\{\big[C(s)\1n+\1n D_1(s)\Th^*_1(s)\1n+\1n
D_2(s)\Th^*_2(s)\big]X(s)\1n+\1n D_1(s)v_1(s)\1n+\1n D_2(s)v^*_2(s)
\1n+\1n\si(s)\Big\}dW(s),\\
\ns\ds\qq\qq\qq\qq\qq\qq\qq\qq\qq\qq\qq\qq\qq\qq s\in[t,T], \\
\ns\ds X(t)=x.\ea\right.\ee
Therefore,
$$\ba{ll}
\ns\ds
J(t,x;\Th^*_1(\cd)X^*(\cd)+v_1^*(\cd),\Th_2^*(\cd)X^*(\cd)+v^*_2(\cd))
\les
J(t,x;\Th^*_1(\cd)X(\cd)+v_1(\cd),\Th^*_2(\cd)X^*(\cd)+v^*_2(\cd))\\
\ns\ds=J(t,x;u_1(\cd),\Th^*_2(\cd)X^*(\cd)+v_2^*(\cd)).\ea$$
Similarly, for any $x\in\dbR^n$, and any $u_2(\cd)\in\cU_2[t,T]$, we can show that
$$J(t,x;\Th^*_1(\cd)X(\cd)+v_1^*(\cd),u_2(\cd))\les
J(t,x;\Th_1^*(\cd)X^*(\cd)+v^*_1(\cd),\Th^*_2(\cd)X^*(\cd)+v_2^*(\cd)).$$
Thus, (i) holds.\endpf

\ms

We note that (iii) of Proposition 3.3 tells us that if we consider
the following state equation (denoting $\Th^*(\cd)=(\Th_1^*(\cd)^T,\Th_2^*(\cd)^T)^T$)
\bel{3.20}\left\{\2n\ba{ll}
\ns\ds dX=\Big\{(A+B\Th^*)X+B_1v_1+B_2v_2^*+b\Big\}ds\\
\ns\ds\qq\qq+\Big\{(C+D\Th^*)X+D_1v_1+D_2v^*_2+\si\Big\}dW(s),\\
\ns\ds X(t)=x,\ea\right.\ee
with the cost functional
\bel{3.21}J_1(t,x;v_1(\cd))=J(t,x;\Th^*_1(\cd)X(\cd)+v_1(\cd),
\Th^*_2(\cd)X(\cd)+v_2^*(\cd)),\ee
then $v_1^*(\cd)$ is an open-loop optimal control of the
corresponding Problem (SLQ). Likewise, if we consider the following
state equation
\bel{3.22}\left\{\2n\ba{ll}
\ns\ds dX=\Big\{(A+B\Th^*)X+B_2v_2+B_1v_1^*+b\Big\}ds\\
\ns\ds\qq\qq+\Big\{(C+D\Th^*)X+D_2v_2+D_1v^*_1+\si\Big\}dW(s),\\
\ns\ds X(t)=x,\ea\right.\ee
with the cost functional
\bel{3.23}J_2(t,x;v_2(\cd))=-J(t,x;\Th^*_1(\cd)X(\cd)+v_1^*(\cd),
\Th^*_2(\cd)X(\cd)+v_2(\cd)),\ee
then $v_2^*(\cd)$ is an open-loop optimal control of the
corresponding Problem (SLQ). The above observation will be useful
below.

\section{Open-Loop Saddle Points and FBSDEs}

In this section, we present a characterization of open-loop saddle
points of Problem (SG) in terms of forward-backward stochastic
differential equations (FBSDE, for short). See \cite{Ma-Yong 1999}
for some relevant results on FBSDEs. The main result of this section
can be stated as follows.

\ms

\bf Theorem 4.1. \sl Let {\rm(SG1)--(SG2)} hold and let
$(t,x)\in[t,T)\times\dbR^n$ be given. Let
$u^*(\cd)\equiv(u_1^*(\cd)^T\2n,u_2^*(\cd)^T\1n)^T\2n\in\cU_1[t,T]\times\cU_2[t,T]$
and $X^*(\cd)\equiv X(\cd\,;t,x,u^*(\cd))$ be the corresponding
state process. Then $u^*(\cd)$ is an open-loop saddle point of
Problem {\rm(SG)} if and only if the following stationarity
conditions hold:
\bel{4.2}B(s)^TY^*(s)+D(s)^TZ^*(s)+S(s)X^*(s)+R(s)u^*(s)+\rho(s)=0,\q
\ae s\in[t,T],~\as\ee
where $(Y^*(\cd),Z^*(\cd))$ is the adapted solution to the following
BSDE:
\bel{BSDE4.1}\left\{\2n\ba{ll}
dY^*(s)=-\big[A(s)^TY^*(s)+C(s)^TZ^*(s)+Q(s)X^*(s)+S(s)^Tu^*(s)+q(s)\big]ds
+Z^*(s)dW(s),\\
\ns\ds\qq\qq\qq\qq\qq\qq\qq\qq\qq\qq\qq\qq\qq\qq\qq\qq\qq s\in[t,T],\\
\ns\ds Y^*(T)=GX^*(T)+g,\ea\right.\ee
and the following convexity-concavity conditions hold: For $i=1,2$,
\bel{convex-concave4.3}\ba{ll}
\ns\ds(-1)^{i-1}\dbE\Big\{\lan GX_i(T),X_i(T)\ran+\int_t^T\Big[\lan
Q(s)X_i(s),X_i(s)\ran+2\lan S_i(s)X_i(s),u_i(s)\ran\\
\ns\ds\qq\qq\qq\qq\qq\qq\qq\qq+\lan
R_{ii}(s)u_i(s),u_i(s)\ran\Big]ds\Big\}\ges0,\qq\forall
u_i(\cd)\in\cU_i[t,T],\ea\ee
where $X_i(\cd)$ solves the following:
\bel{Xi}\left\{\2n\ba{ll}
\ns\ds dX_i(s)=\big[A(s)X_i(s)+B_i(s)u_i(s)\big]ds
+\big[C(s)X_i(s)+D_i(s)u_i(s)\big]dW(s),\q s\in[t,T],\\
\ns\ds X_i(t)=0.\ea\right.\ee

\it Proof. \rm Let
$u^*(\cd)\equiv(u_1^*(\cd)^T,u_2^*(\cd)^T)^T\in\cU_1[t,T]\times\cU_2[t,T]$
and $X^*(\cd)$ be the corresponding state process. Further, let
$(Y^*(\cd),Z^*(\cd))$ be the adapted solution to the BSDE
(\ref{BSDE4.1}). For any $u_1(\cd)\in\cU_1[t,T]$ and $\e\in\dbR$,
let $X^\e(\cd)$ be the solution to the following perturbed state
equation:
$$\left\{\2n\ba{ll}
\ns\ds dX^\e(s)=\big\{A(s)X^\e(s)+B_1(s)[u_1^*(s)+\e
u_1(s)]+B_2(s)u_2^*(s)+b(s)\big\}ds\\
\ns\ds\qq\qq\q+\big\{C(s)X^\e(s)+D_1(s)[u_1^*(s)+\e
u_1(s)]+D_2(s)u_2^*(s)+\si(s)\big\}dW(s),\q s\in[t,T],\\
\ns\ds X^\e(t)=x.\ea\right.$$
Then $X_1(\cd)={X^\e(\cd)-X^*(\cd)\over\e}$ is independent of $\e$
satisfying (\ref{Xi}) (with $i=1$), and
$$\ba{ll}
\ns\ds J(t,x;u^*_1(\cd)+\e u_1(\cd),u_2^*(\cd))
-J(t,x;u_1^*(\cd),u_2^*(\cd))\\
\ns\ds={\e\over2}\dbE\Big\{\lan G\big[2X^*(T)+\e
X_1(T)\big],X_1(T)\ran
+2\lan g,X_1(T)\ran\\
\ns\ds\qq+\int_t^T\Big[\lan\begin{pmatrix}Q&S_1^T&S_2^T\\
                                          S_1&R_{11}&R_{12}\\
                                          S_2&R_{21}&R_{22}\end{pmatrix}
                                   \begin{pmatrix}2X^*+\e X_1\\
                                                  2u_1^*+\e u_1\\
                                                  2u_2^*\end{pmatrix},
                                   \begin{pmatrix}X_1\\
                                                  u_1\\
                                                  0\end{pmatrix}\ran
                                                  +2\lan\begin{pmatrix}q\\
                                                    \rho_1\end{pmatrix},
                                   \begin{pmatrix}X_1\\
                                                  u_1\end{pmatrix}\ran\Big]ds\Big\}\\
\ns\ds=\e\dbE\Big\{\lan GX^*(T)+g,X_1(T)\ran\\
\ns\ds\qq+\int_t^T\Big[\lan QX^*+S^Tu^*+q,X_1\ran+\lan
S_1X^*+R_{11}u_1^*+R_{12}u_2^*+\rho_1,u_1\ran\Big]ds\Big\}\\
\ns\ds\qq+{\e^2\over2}\dbE\Big\{\lan
GX_1(T),X_1(T)\ran+\int_t^T\Big[\lan QX_1,X_1\ran+2\lan
S_1X_1,u_1\ran+\lan R_{11}u_1,u_1\ran\Big]ds\Big\}.\ea$$
On the other hand, we have
$$\ba{ll}
\ns\ds\dbE\Big\{\lan GX^*(T)+g,X_1(T)\ran+\int_t^T\Big[\lan Q
X^*+S^Tu^*+q,X_1\ran+\lan S_1X^*+R_{11}u_1^*+R_{12}u_2^*+\rho_1,u_1\ran\Big]ds\Big\}\\
\ns\ds\q=\dbE\Big\{\int_t^T\Big[\lan-(A^TY^*+C^TZ^*+QX^*+S^Tu^*+q),
X_1\ran+\lan Y^*,AX_1+B_1u_1\ran+\lan Z^*,CX_1+D_1u_1\ran\\
\ns\ds\qq\qq\qq+\lan QX^*+S^Tu^*+q,X_1\ran+\lan
S_1X^*+R_{11}u^*_1+R_{12}u_2^*+\rho_1,u_1\ran\Big]ds\Big\}\\
\ns\ds\q=\dbE\int_t^T\lan
B_1^TY^*+D_1^TZ^*+S_1X^*+R_{11}u_1^*+R_{12}u_2^*+\rho_1,u_1\ran
ds.\ea$$
Hence,
$$\ba{ll}
\ns\ds J(t,x;u^*_1(\cd)+\e u_1(\cd),u_2^*(\cd))
-J(t,x;u_1^*(\cd),u_2^*(\cd))\\
\ns\ds=\e\dbE\Big\{\int_t^T\lan
B_1^TY^*+D_1^TZ^*+S_1X^*+R_{11}u_1^*+R_{12}u_2^*+\rho_1,u_1\ran
ds\Big\}\\
\ns\ds\q+{\e^2\over2}\dbE\Big\{\lan
GX_1(T),X_1(T)\ran+\int_t^T\Big[\lan QX_1,X_1\ran+2\lan
S_1X_1,u_1\ran+\lan R_{11}u_1,u_1\ran\Big]ds\Big\}.\ea$$
Therefore,
$$J(t,x;u_1^*(\cd),u_2^*(\cd))\les J(t,x;u_1^*(\cd)+\e u_1(\cd),u_2^*(\cd)),\qq\forall
u_1(\cd)\in\cU_1[t,T],\q\e\in\dbR,$$
if and only if (\ref{convex-concave4.3}) holds for $i=1$, and
\bel{4.5}B_1^TY^*+D_1^TZ^*+S_1X^*+R_{11}u_1^*+R_{12}u_2^*+\rho_1=0,\q
\ae s\in[t,T],~\as\ee
Similarly,
$$J(t,x;u_1^*(\cd),u_2^*(\cd))\ges J(t,x;u_1^*(\cd),u_2^*(\cd)+\e u_2(\cd)),\qq\forall
u_2(\cd)\in\cU_2[t,T],\q\e\in\dbR,$$
if and only if (\ref{convex-concave4.3}) holds for $i=2$, and
\bel{4.6}B_2^TY^*+D_2^TZ^*+S_2X^*+R_{21}u_1^*+R_{22}u_2^*+\rho_2=0,\q
\ae s\in[t,T],~\as\ee
Combining (\ref{4.5})--(\ref{4.6}), we obtain (\ref{4.2}).\endpf

\ms

From the above result, we see that if Problem (SG) admits an
open-loop saddle point
$u^*(\cd)\equiv(u_1^*(\cd)^T,u_2^*(\cd)^T)^T$, then the following
FBSDE admits an adapted solution
$(X^*(\cd),Y^*(\cd),Z^*(\cd),u^*(\cd))$:
\bel{FBSDE}\left\{\2n\ba{ll}
\ns\ds dX^*(s)\1n=\1n\big[A(s)X^*(s)\1n+\1n B(s)u^*(s)\1n+\1n
b(s)\big]ds\1n+\1n\big[C(s)X^*(s)\1n+\1n
D(s)u^*(s)\1n+\1n\si(s)\big]dW(s),
~s\in[t,T],\\
\ns\ds
dY^*(s)=-\big[A(s)^TY^*(s)+C(s)^TZ^*(s)+Q(s)X^*(s)+S(s)^Tu^*(s)+q(s)\big]ds\\
\ns\ds\qq\qq\qq\qq\qq\qq\qq+Z^*(s)dW(s),\q s\in[t,T],\\
\ns\ds X^*(t)=x,\qq Y^*(T)=GX^*(T)+g,\\
\ns\ds B(s)^TY^*(s)+D(s)^TZ^*(s)+S(s)X^*(s)+R(s)u^*(s)+\rho(s)=0,\q
\ae s\in[t,T],~\as\ea\right.\ee
The following result is concerned with the uniqueness of open-loop
saddle points.

\ms

\bf Theorem 4.2. \sl Let {\rm(SG1)--(SG2)} hold, and let
$(t,x)\in[0,T)\times\dbR^n$ be given. Suppose Problem {\rm(SG)}
admits a unique open-loop saddle point $u^*(\cd)$ at $(t,x)$. Then
FBSDE $(\ref{FBSDE})$ admits a unique adapted solution
$(X^*(\cd),Y^*(\cd),Z^*(\cd),u^*(\cd))$. Conversely, if the
convexity-concavity conditions stated in Theorem 4.1 holds and BFSDE
$(\ref{FBSDE})$ admits a unique adapted solution
$(X^*(\cd),Y^*(\cd),Z^*(\cd),u^*(\cd))$, then $u^*(\cd)$ is the
unique saddle point of Problem {\rm(SG)}.

\ms

\it Proof. \rm Suppose $u^*(\cd)\in\cU[t,T]$ is a unique open-loop
saddle point of Problem (SG). Then by Theorem 4.1, FBSDE
(\ref{FBSDE}) admits an adapted solution
$(X^*(\cd),Y^*(\cd),Z^*(\cd),u^*(\cd))$, and the convex-concave
conditions (\ref{convex-concave4.3}) hold. Now, if FBSDE
(\ref{FBSDE}) admits another different adapted solution $(\bar
X(\cd),\bar Y(\cd),\bar Z(\cd),\bar u(\cd))$. Since the
convexity-concavity conditions are satisfied, by the sufficiency
part of Theorem 4.1, $\bar u(\cd)$ is a different open-loop saddle
point, a contradiction.

\ms

Conversely, if Problem (SG) has two different open-loop saddle
points, then FBSDE (\ref{FBSDE}) will have two different adapted
solutions. \endpf

\section{Closed-Loop Saddle Points and Riccati Equations}

We now look at closed-loop saddle points for Problem (SG). First, we
present the following result which is a consequence of Theorem 4.1.

\ms

\bf Proposition 5.1. \sl Let {\rm(SG1)--(SG2)} hold. Let
$(\Th^*(\cd),v^*(\cd))\in\sQ[t,T]\times\cU[t,T]$ be a closed-loop
saddle point of Problem {\rm(SG)}. Then the following FBSDE admits
an adapted solution $(X^*(\cd),Y^*(\cd),Z^*(\cd))$:
\bel{FBSDE5.1}\left\{\2n\ba{ll}
\ns\ds
dX^*=\big\{(A+B\Th^*)X^*+Bv^*+b\big\}ds+\big\{(C+D\Th^*)X^*+Dv^*+\si\big\}dW(s),
\q s\in[t,T],\\
\ns\ds
dY^*(s)=-\big\{A^TY^*+C^TZ^*+(Q+S^T\Th^*)X^*+S^Tv^*+q\big\}ds+Z^*dW(s),\\
\ns\ds X^*(t)=x,\qq Y^*(T)=GX^*(T)+g,\ea\right.\ee
and the following stationarity condition holds:
\bel{}Rv^*+B^TY^*+D^TZ^*+(S+R\Th^*)X^*+\rho=0,\q \ae~\as\ee

\it Proof. \rm Let $(\Th^*(\cd),v^*(\cd))\in\sQ[t,T]\times\cU[t,T]$
be a closed-loop saddle point of Problem (SG) with
$\Th^*(\cd)=(\Th_1^*(\cd)^T,\Th_2^*(\cd)^T)^T$ and
$v^*(\cd)=(v_1^*(\cd)^T,v_2^*(\cd)^T)^T$. We consider state equation
(\ref{3.20}) with the cost functional (\ref{3.21}) for which we
carry out some computation: (denoting $\wt v=(v_1^T,(v_2^*)^T)^T$)
$$\ba{ll}
\ns\ds J_1(t,x;v_1(\cd))\equiv
J(t,x;\Th^*X(\cd)+\wt v(\cd))\\
\ns\ds={1\over2}\dbE\Big\{\lan GX(T),X(T)\ran+2\lan
g,X(T)\ran+\int_t^T\Big[\lan QX,X\ran+2\lan
SX,\Th^*X+\wt v\ran\\
\ns\ds\qq\q+\lan R(\Th^*X+\wt v),\Th^*X+\wt v\ran+2\lan
q,X\ran+2\lan\rho,\Th^*X+\wt v\ran\Big]
ds\Big\}\\
\ns\ds={1\over2}\dbE\Big\{\lan GX(T),X(T)\ran+2\lan
g,X(T)\ran+\int_t^T\Big[\lan[Q+(\Th^*)^TS+S^T\Th^*+(\Th^*)^TR\Th^*]X,X\ran\\
\ns\ds\qq\q+2\lan\begin{pmatrix}(S_1+R_1\Th^*)X\\ (S_2+R_2\Th^*)X\end{pmatrix},\begin{pmatrix}v_1\\
v_2^*\end{pmatrix}\ran+\lan\begin{pmatrix}R_{11}&R_{12}\\
R_{21}&R_{22}\end{pmatrix}\begin{pmatrix}v_1\\
v_2^*\end{pmatrix},\begin{pmatrix}v_1\\
v_2^*\end{pmatrix}\ran\\
\ns\ds\qq\q+2\lan
q+(\Th^*)^T\rho,X\ran+2\lan\begin{pmatrix}\rho_1\\ \rho_2\end{pmatrix},\begin{pmatrix}v_1\\
v_2^*\end{pmatrix}\ran\Big] ds\Big\}\\
\ns\ds={1\over2}\dbE\Big\{\lan GX(T),X(T)\ran+2\lan
g,X(T)\ran+\int_t^T\Big[\lan[Q+(\Th^*)^TS+S^T\Th^*+(\Th^*)^TR\Th^*]X,X\ran\\
\ns\ds\qq\q+2\lan(S_1+R_1\Th^*)X,v_1\ran+2\lan
q+(\Th^*)^T\rho+(S_2+R_2\Th^*)^Tv_2^*,X\ran\\
\ns\ds\qq\q+\lan R_{11}v_1,v_1\ran+2\lan\rho_1+R_{12}v_2^*,v_1\ran
+\lan R_{22}v_2^*,v_2^*\ran+2\lan\rho_2,v_2^*\ran\Big]
ds\Big\}.\ea$$
We know that $v_1^*(\cd)$ is an open-loop optimal control for the
problem with state equation (\ref{3.20}) and the above cost
functional. Thus, according to Theorem 4.1, we have
$$\ba{ll}
\ns\ds0=B^T_1Y^*+D_1^TZ^*+(S_1+R_1\Th^*)X^*+R_{11}v_1^*+\rho_1+R_{12}v_2^*,\q \ae~\as\ea$$
with $(Y^*(\cd),Z^*(\cd))$ being the adapted solution to the
following BSDE:
$$\ba{ll}
\ns\ds
dY^*=-\big\{(A+B\Th^*)^TY^*+(C+D\Th^*)^TZ^*+[Q+(\Th^*)^TS+S^T\Th^*+(\Th^*)^TR\Th^*]
X^*\\
\ns\ds\qq\qq+(S_1+R_1\Th^*)^Tv_1^*+q+(\Th^*)^T\rho+(S_2+R_2\Th^*)^Tv_2^*\big\}ds+Z^*dW\\
\ns\ds\qq=-\big\{A^TY^*+C^TZ^*+QX^*+S^T(\Th^*X^*+v^*)+q\\
\ns\ds\qq\qq+(\Th^*)^T[B^TY^*+D^TZ^*+SX^*+R(\Th^*X^*+v^*)+\rho]\big\}ds+Z^*dW\\
\ns\ds Y^*(T)=GX^*(T)+g.\ea$$
Likewise, by considering state equation (\ref{3.22}) and payoff
functional (\ref{3.23}), we can obtain
$$\ba{ll}
\ns\ds0=B^T_2Y^*+D_2^TZ^*+(S_2+R_2\Th^*)X^*+R_{21}v_1^*+\rho_2+R_{22}v_2^*,\q \ae~\as\ea$$
with $(Y^*(\cd),Z^*(\cd))$ being the adapted solution to the same
BSDE as above. Thus,
$$\ba{ll}
\ns\ds0=B^TY^*+D^TZ^*+(S+R\Th^*)X^*+Rv^*+\rho,\q \ae~\as\ea$$
Then the above BSDE is reduced to that in (\ref{FBSDE5.1}). \endpf

\ms

The following result gives a characterization for closed-loop saddle
points of Problem (SG).

\ms

\bf Theorem 5.2. \sl Let {\rm(SG1)--(SG2)} hold. Then Problem
{\rm(SG)} admits a closed-loop saddle point
$(\Th^*(\cd),v^*(\cd))\in\sQ[t,T]\times\cU[t,T]$ with
$\Th^*(\cd)\equiv(\Th_1^*(\cd)^T,\Th_2^*(\cd)^T)^T$ and
$v^*(\cd)\equiv(v_1^*(\cd)^T,v_2^*(\cd)^T)^T$ if and only if the
following Riccati equation:
\bel{Riccati}\left\{\2n\ba{ll}
\ns\ds\dot P(s)+P(s)A(s)+A(s)^TP(s)+C(s)^TP(s)C(s)+Q(s)\\
\ns\ds\q-\big[P(s)B(s)+C(s)^TP(s)D(s)+S(s)^T\big]\big[R(s)+D(s)^T
P(s)D(s)\big]^\dag\\
\ns\ds\qq\cd\big[B(s)^TP(s)
+D(s)^TP(s)C(s)+S(s)\big]=0,\qq \ae s\in[t,T],\\
\ns\ds P(T)=G,\ea\right.\ee
admits a solution $P(\cd)\in C([t,T];\dbS^n)$ such that
\bel{5.2}\cR\big(B(s)^TP(s)+D(s)^TP(s)C(s)+S(s)\big)\subseteq\cR\big(R(s)+D(s)^TP(s)D(s)\big),\q
\ae s\in[t,T],\ee
\bel{th
L2}\big[R(\cd)+D(\cd)^TP(\cd)D(\cd)\big]^\dag\big[B(\cd)^TP(\cd)
+D(\cd)^TP(\cd)C(\cd)+S(\cd)\big] \in L^2(t,T;\dbR^{m\times n}),\ee
\bel{5.6}R_{11}(s)+D_1(s)^TP(s)D_1(s)\ges0,\q
R_{22}(s)+D_2(s)^TP(s)D_2(s)\les0,\q\ae s\in[t,T],\ee
and the adapted solution $(\eta(\cd),\z(\cd))$ of the
following BSDE:
\bel{5.8}\left\{\2n\ba{ll}
\ns\ds d\eta=-\Big\{\big[A^T-(PB+C^TPD+S^T)(R+D^TPD)^\dag
B^T\big]\eta\\
\ns\ds\qq\qq+\big[C^T-(PB+C^TPD+S^T)(R+D^TPD)^\dag D^T\big]\z\\
\ns\ds\qq\qq+\big[C^T-(PB+C^TPD+S^T)(R+D^TPD)^\dag
D^T\big]P\si\\
\ns\ds\qq\qq-(PB+C^TPD+S^T)(R+D^TPD)^\dag
\rho+Pb+q\Big\}ds+\z dW(s),\\
\ns\ds\eta(T)=g,\ea\right.\ee
satisfies
\bel{}B(s)^T\eta(s)+D(s)^T\z(s)+D(s)^TP(s)\si(s)+\rho(s)\in\cR\big(R(s)+D(s)^TP(s)D(s)\big),
\q \ae s\in[t,T],~\as\ee
and
\bel{v
L2}\big[R(\cd)+D(\cd)^TP(\cd)D(\cd)\big]^\dag\big[B(\cd)^T\eta(\cd)+D(\cd)^T\z(\cd)
+D(\cd)^TP(\cd)\si(\cd)+\rho(\cd)\big]\in L_\dbF^2(t,T;\dbR^m).\ee
In this case, the closed-loop saddle point $(\Th^*(\cd),v^*(\cd))$
admits the following representation:
\bel{5.10}\left\{\2n\ba{ll}
\ns\ds\Th^*(\cd)=-[R(\cd)+D(\cd)^TP(\cd)D(\cd)]^\dag[B(\cd)^TP(\cd)
+D(\cd)^TP(\cd)C(\cd)+S(\cd)]\\
\ns\ds\qq\qq+\big\{I-\big[R(\cd)+D(\cd)^TP(\cd)D(\cd)\big]^\dag
\big[R(\cd)+D(\cd)^TP(\cd)D(\cd)\big]\big\}\th(\cd),\\
\ns\ds
v^*(\cd)=-\big[R(\cd)+D(\cd)^TP(\cd)D(\cd)\big]^\dag\big[B(\cd)^T\eta(\cd)
+D(\cd)^T\z(\cd)+D(\cd)^TP(\cd)\si(\cd)+\rho(\cd)\big]\\
\ns\ds\qq\qq+\big\{I-\big[R(\cd)+D(\cd)^TP(\cd)D(\cd)\big]^\dag
\big[R(\cd)+D(\cd)^TP(\cd)D(\cd)]\big\}\n(\cd),\ea\right.\ee
for some $\th(\cd)\in L^2(t,T;\dbR^{m\times n})$ and $\n(\cd)\in L_\dbF^2(t,T;\dbR^m)$.

\ms

Further, the value function admits the following representation:
\bel{Value}\ba{ll}
\ns\ds V(t,x)={1\over2}\,\dbE\Big\{\lan
P(t)x,x\ran+2\lan\eta(t),x\ran+\1n\int_t^T\3n\big[\lan
P\si,\si\ran\1n+\1n2\lan\eta,b\ran+\1n2\lan\z,\si\ran\\
\ns\ds\qq\qq\qq\q-\1n \lan(R\1n+\1n D^T\1n
PD)^\dag(B^T\1n\eta\1n+D^T\1n\z\1n+\1n D^T\1n P\si\1n+\1n\rho),B^T\1n\eta\1n+D^T\1n\z\1n+\1n
D^T\1n P\si\1n+\1n\rho\ran\big]ds\Big\}.\ea\ee

\ms

\it Proof. \rm {\it Necessity.} Let $(\Th^*(\cd),v^*(\cd))$ be a
closed-loop saddle point of Problem (SG) over $[t,T]$, where
$\Th^*(\cd)\equiv(\Th^*_1(\cd)^T,\Th^*_2(\cd)^T)^T\in\sQ_1[t,T]\times\sQ_2[t,T]$
and
$v^*(\cd)\equiv(v_1^*(\cd)^T,v_2^*(\cd)^T)^T\in\cU_1[t,T]\times\cU_2[t,T]$.
Then, by Proposition 5.1, for any $x\in\dbR^n$, the following FBSDE
admits an adapted solution $(X^*(\cd),Y^*(\cd),Z^*(\cd))$:
\bel{}\left\{\2n\ba{ll}
\ns\ds
dX^*=\big\{(A+B\Th^*)X^*+Bv^*+b\big\}ds+\big\{(C+D\Th^*)X^*+Dv^*+\si\big\}dW(s),
\q s\in[t,T],\\
\ns\ds
dY^*(s)=-\big\{A^TY^*+C^TZ^*+(Q+S^T\Th^*)X^*+S^Tv^*+q\big\}ds+Z^*dW(s),\\
\ns\ds X^*(t)=x,\qq Y^*(T)=GX^*(T)+g,\ea\right.\ee
and the following stationarity condition holds:
\bel{5.9}B^TY^*+D^TZ^*+(S+R\Th^*)X^*+Rv^*+\rho=0,\q\ae~\as\ee
Since the above admits a solution for each $x\in\dbR^n$, and
$(\Th^*(\cd),v^*(\cd))$ is independent of $x$, by subtracting
solutions corresponding $x$ and $0$, the later from the former, we
see that for any $x\in\dbR^n$, as long as $(X(\cd),Y(\cd),Z(\cd))$
is the adapted solution to the following FBSDE:
\bel{FBSDE-closed}\left\{\2n\ba{ll}
\ns\ds dX=(A+B\Th^*)Xds+(C+D\Th^*)XdW(s),\q s\in[t,T],\\
\ns\ds
dY=-\big[A^TY+C^TZ+(Q+S^T\Th^*)X\big]ds+ZdW(s),\q s\in[t,T],\\
\ns\ds X(t)=x,\qq Y(T)=GX(T),\ea\right.\ee
one must have the following stationarity condition:
\bel{5.11}B^TY+C^TZ+(S+R\Th^*)X=0,\q \ae s\in[t,T],~\as\ee
Now, we let
\bel{}\left\{\2n\ba{ll}
\ns\ds d\dbX(s)=\big[A(s)+B(s)\Th^*(s)\big]\dbX(s)ds
+\big[C(s)+D(s)\Th^*(s)\big]\dbX(s)dW(s),\q s\in[t,T],\\
\ns\ds\dbX(t)=I,\ea\right.\ee
and let
\bel{}\left\{\2n\ba{ll}
\ns\ds
d\dbY(s)=\big\{-A(s)^T\dbY(s)-C(s)^T\dbZ(s)-\big[Q(s)+S(s)^T\Th^*(s)\big]
\dbX(s)\big\}ds\\
\ns\ds\qq\qq\qq\qq\qq\qq\qq\qq\qq+\dbZ(s)dW(s),\qq s\in[t,T],\\
\ns\ds\dbY(T)=G\dbX(T).\ea\right.\ee
Clearly, $\dbX(\cd)$, $\dbY(\cd)$, and $\dbZ(\cd)$ are all
well-defined square matrix valued processes. Further, (\ref{5.11})
implies
\bel{5.14}B(s)^T\dbY(s)+D(s)^T\dbZ(s)+\big[S(s)+R(s)\Th^*(s)\big]\dbX(s)=0,
\q \ae s\in[t,T],~\as\ee
Clearly, $\dbX(\cd)^{-1}$ exists and satisfies the following:
\bel{}\left\{\2n\ba{ll}
\ns\ds
d\big[\dbX(s)^{-1}]=\dbX(s)^{-1}\big\{\big[C(s)+D(s)\Th^*(s)\big]^2
-A(s)-B(s)\Th^*(s)\big\}ds\\
\ns\ds\qq\qq\qq-\dbX(s)^{-1}\big[C(s)+D(s)\Th^*(s)\big]dW(s),\qq s\in[t,T],\\
\ns\ds\dbX(t)^{-1}=I.\ea\right.\ee
We define
$$P(\cd)=\dbY(\cd)\dbX(\cd)^{-1},\qq\Pi(\cd)=\dbZ(\cd)\dbX(\cd)^{-1}.$$
Then (\ref{5.14}) implies
\bel{5.16}B^TP+D^T\Pi+(S+R\Th^*)=0,\q\ae\ee
Also, by It\^o's formula,
$$\ba{ll}
\ns\ds
dP\1n=\1n\Big\{\1n-\1n\big[A^T\dbY+C^T\dbZ+(Q+S^T\Th^*)\dbX\big]\dbX^{-1}\2n
+\1n\dbY\dbX^{-1}\big[(C\1n+D\Th^*)^2-A-B\Th^*\big]
-\dbZ\dbX^{-1}(C\1n+D\Th^*)\Big\}ds\\
\ns\ds\qq\qq+\Big\{\dbZ\dbX^{-1}-\dbY\dbX^{-1}(C+D\Th^*)\Big\}dW(s)\\
\ns\ds\qq=\Big\{-A^TP-C^T\Pi-Q-S^T\Th^*
+P\big[(C\1n+D\Th^*)^2-A-B\Th^*\big]
-\Pi(C\1n+D\Th^*)\Big\}ds\\
\ns\ds\qq\qq+\Big\{\Pi-P(C+D\Th^*)\Big\}dW(s).\ea$$
Let
$$\L=\Pi-P(C+D\Th^*).$$
Then
$$\ba{ll}
\ns\ds dP=\Big\{-A^TP-C^T[\L+P(C+D\Th^*)]-Q-S^T\Th^*
+P\big[(C\1n+D\Th^*)^2-A-B\Th^*\big]\\
\ns\ds\qq\qq-[\L+P(C+D\Th^*)](C\1n+D\Th^*)\Big\}ds+\L
dW(s)\\
\ns\ds\q~=\1n\Big\{\1n-\1n PA\1n-\1n A^T\1n P\1n-\1n\L C\1n-\1n
C^T\1n\L\1n-\1n C^T\1n PC\1n-\1n(PB\1n+\1n C^T\1n PD\1n+\1n
S^T\2n+\1n\L D)\Th^*\1n-\1n Q\Big\}ds\1n+\1n\L dW(s),\ea$$
and $P(T)=G$. Thus, $(P(\cd),\L(\cd))$ is the adapted solution of a
BSDE with deterministic coefficients. Hence, $P(\cd)$ is
deterministic and $\L(\cd)=0$ which means
\bel{}\Pi=\dbZ\dbX^{-1}=P(C+D\Th^*).\ee
Therefore,
\bel{2.65}\dot P+PA+A^TP+C^TPC+(PB+C^TPD+S^T)\Th^*+Q=0,\q\ae\ee
and (\ref{5.16}) becomes
\bel{5.20}0=B^TP+D^TP(C+D\Th^*)+S+R\Th^*
=B^TP+D^TPC+S+(R+D^TPD)\Th^*,\q\ae\ee
This implies
$$\cR\big(B^TP+D^TPC+S\big)\subseteq\cR\big(R+D^TPD\big),\q\ae$$
Using (\ref{5.20}), (\ref{2.65}) can be written as
$$\ba{ll}
\ns\ds0=\dot P+P(A+B\Th^*)+(A+B\Th^*)^TP+(C+D\Th^*)^TP(C+D\Th^*)\\
\ns\ds\qq\q+(\Th^*)^TR\Th^*+S^T\Th^*+(\Th^*)^TS+Q,\q\ae\ea$$
Since $P(T)=G\in\dbS^n$ and $Q(\cd),R(\cd)$ are symmetric, by uniqueness,
we must have $P(\cd)\in C([t,T];\dbS^n)$.
Denoting $\h R=R+D^TPD$, since
$$\h R^\dag(B^TP+D^TPC+S)=-\h R^\dag\h R\Th^*,$$
and $\h R^\dag\h R$ is an orthogonal projection,
we see that (\ref{th L2}) holds and
$$\Th^*=-\h R^\dag(B^TP+D^TPC+S)+\big(I-\h R^\dag\h R\,\big)\th$$
for some $\th(\cd)\in L^2(t,T;\dbR^{m\times n})$. Consequently,
\bel{5.21}\ba{ll}
\ns\ds(PB+C^TPD+S^T)\Th^*=(\Th^*)^T\h R\h
R^\dag(B^TP+D^TPC+S)\\
\ns\ds=-(PB+C^TPD+S^T)\h R^\dag(B^TP+D^TPC+S).\ea\ee
Plug the above into (\ref{2.65}), we obtain Riccati equation
(\ref{Riccati}). To determine $v^*(\cd)$, we define
$$\left\{\2n\ba{ll}
\ns\ds\eta=Y^*-PX^*,\\
\ns\ds\z=Z^*-P(C+D\Th^*)X^*-PDv^*-P\si.\ea\right.\qq s\in[t,T].$$
Then
$$\ba{ll}
\ns\ds d\eta=dY^*-\dot PX^*ds-PdX^*\\
\ns\ds\qq=-\big[A^TY^*+C^TZ^*+(Q+S^T\Th^*)X^*+S^Tv^*+q\big]ds+Z^*dW\\
\ns\ds\qq\q+\Big\{\big[PA\1n+\1n A^T\1n P\1n+\1n C^T\1n PC\1n+\1n
Q\1n-\1n(PB\1n+C^T\1n PD\1n+\1n S^T)(R\1n+\1n D^T\1n PD)^\dag(B^T\1n
P\1n+\1n D^T\1n PC\1n+\1n S)\big]X^*\\
\ns\ds\qq\q-P[(A+B\Th^*)X^*+Bv^*+b]\Big\}ds-P\big[(C+D\Th^*)X^*+Dv^*+\si]dW\\
\ns\ds\qq=-\Big\{A^T(\eta+PX^*)+C^T\big[\z+P(C+D\Th^*)X^*+PDv^*+P\si\big]
+(Q+S^T\Th^*)X^*+S^Tv^*+q\\
\ns\ds\qq\q-\big[PA\1n+\1n A^T\1n P\1n+\1n C^T\1n PC\1n+\1n
Q\1n-\1n(PB\1n+\1n C^T\1n PD\1n+\1n S^T)(R\1n+\1n D^T\1n
PD)^\dag(B^T\1n P\1n+\1n D^T\1n PC\1n+\1n S)\big]X^*\\
\ns\ds\qq\q+P[(A+B\Th^*)X^*+Bv^*+b]\Big\}ds+\z dW\\
\ns\ds\qq=\Big\{\1n-\1n A^T\1n\eta-C^T\1n\z\1n-\1n(PB\1n+\1n C^T\1n
PD\1n+\1n S^T)\Th^*X^*\1n-\1n(PB\1n+\1n C^T\1n PD\1n+\1n
S^T)v^*\1n-\1n C^T\1n P\si\1n
-\1n Pb\1n-\1n q\\
\ns\ds\qq\q-\big[(PB+C^TPD+S^T)(R+D^TPD)^\dag(B^TP+D^TPC+S)\big]X^*\Big\}ds+\z dW\\
\ns\ds\qq=-\big[A^T\eta+C^T\z+(PB+C^TPD+S^T)v^*+C^TP\si+Pb+q\big]ds+\z
dW.\ea$$
According to (\ref{5.9}), we have
$$\ba{ll}
\ns\ds0=B^TY^*+D^TZ^*+(S+R\Th^*)X^*+Rv^*+\rho\\
\ns\ds\q=B^T(\eta+PX^*)+D^T[\z+P(C+D\Th^*)X^*+PDv^*+P\si]+(S+R\Th^*)X^*
+Rv^*+\rho\\
\ns\ds\q=[B^TP+D^TPC+S+(R+D^TPD)\Th^*]X^*+B^T\eta+D^T\z+D^TP\si+\rho
+(R+D^TPD)v^*\\
\ns\ds\q=B^T\eta+D^T\z+D^TP\si+\rho+(R+D^TPD)v^*.\ea$$
Hence,
$$B^T\eta+D^T\z+D^TP\si+\rho\in\cR(R+D^TPD),\q\ae~\as$$
Since
$\h R^\dag(B^T\eta+D^T\z+D^TP\si+\rho)=-\h R^\dag\h Rv^*,$
and $\h R^\dag\h R$ is an orthogonal projection,
we see that (\ref{v L2}) holds and
$$v^*=-\h R^\dag(B^T\eta+D^T\z+D^TP\si+\rho)+\big[I-\h R^\dag\h R\big]
\n$$
for some $\n(\cd)\in L_\dbF^2(t,T;\dbR^m)$. Consequently,
$$\ba{ll}
\ns\ds(PB+C^TPD+S^T)v^*=-(PB+C^TPD+S^T)(R+D^TPD)^\dag(B^T\eta+D^T\z+D^TP\si+\rho)\\
\ns\ds\qq\qq\qq\qq\qq\qq+(PB+C^TPD+S^T)\big[I-(R+D^TPD)^\dag(R+D^TPD)\big]\n\\
\ns\ds\qq\qq\qq\qq\qq=-(PB+C^TPD+S^T)(R+D^TPD)^\dag(B^T\eta+D^T\z+D^TP\si+\rho).\ea$$
Then
$$\ba{ll}
\ns\ds A^T\eta+C^T\z+(PB+C^TPD+S^T)v^*+C^TP\si+Pb+q\\
\ns\ds=A^T\1n\eta\1n+\1n C^T\1n\z\1n-\1n(PB\1n+\1n C^T\1n PD\1n+\1n
S^T) (R\1n+\1n D^T\1n PD)^\dag(B^T\1n\eta\1n+\1n D^T\1n\z\1n+\1n
D^T\1n P\si
\1n+\1n\rho)\1n+\1n C^T\1n P\si\1n+\1n Pb\1n+\1n q\\
\ns\ds=\1n\big[A^T\2n-\1n(PB\1n+\1n C^T\1n PD+\1n S^T)(R\1n+\1n
D^TPD)^\dag B^T\big]\eta\1n+\1n\big[C^T\2n-(PB\1n+C^T\1n
PD\1n+\1n S^T\1n)(R\1n+\1n D^T\1n PD)^\dag D^T\big]\z\\
\ns\ds\qq+\big[C^T-(PB+C^TPD+S^T)(R+D^TPD)^\dag
D^T\big]P\si\\
\ns\ds\qq-(PB+C^TPD+S^T)(R+D^TPD)^\dag\rho+Pb+q.\ea$$
Therefore, $(\eta,\z)$ is the adapted solution to the following BSDE:
$$\left\{\2n\ba{ll}
\ns\ds d\eta=-\Big\{\big[A^T-(PB+C^TPD+S^T)(R+D^TPD)^\dag
B^T\big]\eta\\
\ns\ds\qq\qq+\big[C^T-(PB+C^TPD+S^T)(R+D^TPD)^\dag D^T\big]\z\\
\ns\ds\qq\qq+\big[C^T-(PB+C^TPD+S^T)(R+D^TPD)^\dag
D^T\big]P\si\\
\ns\ds\qq\qq-(PB+C^TPD+S^T)(R+D^TPD)^\dag
\rho+Pb+q\Big\}ds+\z dW(s),\\
\ns\ds\eta(T)=g.\ea\right.$$
To prove (\ref{5.6}) as well as the sufficiency, we take any
$u(\cd)=(u_1(\cd)^T,u_2(\cd)^T)^T\in\cU_1[t,T]\times\cU_2[t,T]$, let
$X(\cd)\equiv X(\cd\,;t,x,u(\cd))$ be the corresponding state
process. Then
$$\ba{ll}
\ns\ds J(t,x;u(\cd))={1\over2}\dbE\Big\{\lan GX(T),X(T)\ran+2\lan
g,X(T)\ran\\
\ns\ds\qq\qq\qq\qq+\int_t^T\big[\lan QX,X\ran+2\lan SX,u\ran+\lan
Ru,u\ran+2\lan q,X\ran+2\lan\rho,u\ran\big]ds\Big\}\ea$$
$$\ba{ll}
\ns\ds={1\over2}\dbE\Big\{\lan P(t)x,x\ran+2\lan\eta(t),x\ran\\
\ns\ds\qq+\2n\int_t^T\3n\big\{\1n\lan\big[-\1n PA\1n-\1n A^T\1n
P\1n-\1n C^T\1n PC\1n-\1n Q\1n+\1n(PB\1n+\1n C^T\1n PD\1n+\1n
S^T)(R\1n+\1n D^T\1n PD)^\dag(B^T\1n P\1n
+\1n D^T\1n PC\1n+\1n S)\big]X,X\ran\\
\ns\ds\qq+\lan P(AX\1n+\1n Bu\1n+\1n b),X\ran\1n+\1n\lan
PX,AX\1n+\1n Bu\1n+\1n b\ran\1n+\1n\lan P(CX\1n+\1n
Du\1n+\1n\si),CX\1n+\1n Du
\1n+\1n\si\ran\\
\ns\ds\qq+2\lan\big[-A^T+(PB+C^TPD+S^T)(R+D^TPD)^\dag
B^T]\eta,X\ran\\
\ns\ds\qq+2\lan\big[-C^T+(PB+C^TPD+S^T)(R+D^TPD)^\dag
D^T]P\si,X\ran\\
\ns\ds\qq+2\lan\big[-C^T+(PB+C^TPD+S^T)(R+D^TPD)^\dag
D^T]\z,X\ran\\
\ns\ds\qq+2\lan\big(PB+C^TPD+S^T)(R+D^TPD)^\dag
\rho-Pb-q,X\ran+2\lan\z,CX+Du+\si\ran\\
\ns\ds\qq+2\lan\eta,AX+Bu+b\ran+\lan QX,X\ran+2\lan SX,u\ran+\lan
Ru,u\ran+2\lan q,X\ran+2\lan\rho,u\ran\big\}ds\Big\}\\
%
%
%
%
%
%
%
\ns\ds={1\over2}\dbE\Big\{\lan
P(t)x,x\ran+2\lan\eta(t),x\ran+\int_t^T\Big[
\lan P\si,\si\ran+2\lan\eta,b\ran+2\lan\z,\si\ran\\
\ns\ds\qq+\lan(PB\1n+\1n C^T\1n PD\1n+\1n S^T)(R\1n+\1n D^T\1n
PD)^\dag(B^T\1n P\1n
+\1n D^T\1n PC\1n+\1n S)X,X\ran\\
\ns\ds\qq+2\lan(B^TP+D^TPC+S)X+B^T\eta+D^T\z+D^TP\si+\rho,u\ran\1n
+\1n\lan(R+D^TPD)u,u\ran\\
\ns\ds\qq+2\lan(PB+C^TPD+S^T)(R+D^TPD)^\dag
(B^T\eta+D^T\z+D^TP\si+\rho),X\ran\Big]ds\Big\}.\ea$$
Note that
$$\left\{\2n\ba{ll}
\ns\ds B^TP+D^TPC+S=-(R+D^TPD)\Th^*\equiv-\h R\Th^*,\\
\ns\ds B^T\eta+D^T\z+D^TP\si+\rho=-(R+D^TPD)v^*\equiv-\h Rv^*.\ea\right.$$
Also, one has
$$\ba{ll}
\ns\ds\lan(R+D^TPD)v^*,v^*\ran=\lan\h R\h
R^\dag(B^T\eta+D^T\z+D^TP\si+\rho),\h
R^\dag(B^T\eta+D^T\z+D^TP\si+\rho)\ran\\
\ns\ds\qq\qq\qq\qq\q~=\lan(R+D^TPD)^\dag(B^T\eta+D^T\z+D^TP\si+\rho),B^T\eta+D^T\z+D^TP\si+\rho\ran.\ea$$
Thus,
$$\ba{ll}
\ns\ds J(t,x;u(\cd))={1\over2}\dbE\Big\{\lan
P(t)x,x\ran+2\lan\eta(t),x\ran+\int_t^T\Big[
\lan P\si,\si\ran+2\lan\eta,b\ran+2\lan\z,\si\ran\\
\ns\ds\qq+\lan(PB\1n+\1n C^T\1n PD\1n+\1n S^T)(R\1n+\1n D^T\1n
PD)^\dag(B^T\1n P\1n
+\1n D^T\1n PC\1n+\1n S)X,X\ran\\
\ns\ds\qq+2\lan(B^TP+D^TPC+S)X+B^T\eta+D^T\z+D^TP\si+\rho,u\ran\1n
+\1n\lan(R+D^TPD)u,u\ran\\
\ns\ds\qq+2\lan(PB+C^TPD+S^T)(R+D^TPD)^\dag
(B^T\eta+D^T\z+D^TP\si+\rho),X\ran\Big]ds\Big\}\\
\ns\ds={1\over2}\dbE\Big\{\lan
P(t)x,x\ran+2\lan\eta(t),x\ran+\int_t^T\Big[
\lan P\si,\si\ran+2\lan\eta,b\ran+2\lan\z,\si\ran\\
\ns\ds\qq+\lan(\Th^*)^T\h R\h R^\dag\h R\Th^* X,X\ran-2\lan\h
R(\Th^* X+v^*),u\ran\1n +\1n\lan\h Ru,u\ran+2\lan(\Th^*)^T\h
R\h R^\dag\h
Rv^*,X\ran\Big]ds\Big\}\\
\ns\ds={1\over2}\dbE\Big\{\lan
P(t)x,x\ran+2\lan\eta(t),x\ran+\int_t^T\Big[
\lan P\si,\si\ran+2\lan\eta,b\ran+2\lan\z,\si\ran\\
\ns\ds\qq+\lan\h R\Th^*X,\Th^*X\ran-2\lan\h
R(\Th^*X+v^*),u\ran\1n+\1n\lan\h Ru,u\ran+2\lan\h R\Th^*
X,v^*\ran\Big]ds\Big\}\\
\ns\ds={1\over2}\dbE\Big\{\lan
P(t)x,x\ran+2\lan\eta(t),x\ran+\int_t^T\Big[\lan P\si,\si\ran+2\lan\eta,b\ran+2\lan\z,\si\ran\\
\ns\ds\qq\qq-
\lan(R+D^TPD)^\dag(B^T\eta+D^T\z+D^TP\si+\rho),B^T\eta+D^T\z+D^TP\si+\rho\ran\\
\ns\ds\qq\qq+\lan(R+D^TPD)(u-\Th^*X-v^*),u-\Th^*X-v^*\ran\Big]ds\Big\}\\
\ns\ds=J\big(t,x;\Th^*(\cd)X^*(\cd)+v^*(\cd)\big)
+{1\over2}\dbE\int_t^T\lan(R+D^TPD)(u-\Th^*X-v^*),u-\Th^*X-v^*\ran
ds.\ea$$
Consequently,
$$\ba{ll}
\ns\ds
J(t,x;\Th_1^*(\cd)X(\cd)+v_1(\cd),\Th_2^*(\cd)X(\cd)+v_2^*(\cd))\\
\ns\ds=J(t,x;\Th^*(\cd)X^*(\cd)+v^*(\cd))
+{1\over2}\dbE\int_t^T\lan(R_{11}+D_1^TPD_1)(v_1-v_1^*),v_1-v_1^*\ran
ds.\ea$$
Hence,
$$J(t,x;\Th^*(\cd)X^*(\cd)+v^*(\cd))\les J(t,x;\Th_1^*(\cd)X(\cd)+v_1(\cd),
\Th^*_2(\cd)X(\cd)+v_2^*(\cd)),\q\forall v_1(\cd)\in\cU_1[t,T],$$
if and only if
$$R_{11}+D_1^TPD_1\ges0,\q\ae s\in[t,T].$$
Similarly,
$$\ba{ll}
\ns\ds
J(t,x;\Th_1^*(\cd)X(\cd)+v_1^*(\cd),\Th_2^*(\cd)X(\cd)+v_2(\cd))\\
\ns\ds=J(t,x;\Th^*(\cd)X^*(\cd)+v^*(\cd))
+{1\over2}\dbE\int_t^T\lan(R_{22}+D_2^TPD_2)(v_2-v_2^*),v_2-v_2^*\ran
ds.\ea$$
Hence,
$$J(t,x;\Th^*(\cd)X^*(\cd)+v^*(\cd))\ges J(t,x;\Th_1^*(\cd)X(\cd)+v^*_1(\cd),
\Th^*_2(\cd)X(\cd)+v_2(\cd)),\q\forall v_2(\cd)\in\cU_2[t,T],$$
if and only if
$$R_{22}+D_2^TPD_2\les0,\q\ae s\in[t,T].$$
That is, $(\Th^*(\cd),v^*(\cd))$ is a closed-loop saddle point of
Problem (SG). \endpf

\ms

Note that by letting $m_2=0$, we obtain a proof for Theorem 2.4.

\ms

We point out here that the solution of the Riccati equation
(\ref{Riccati}) may be non-unique. We will present an example in
Section 6. A solution $P(\cd)$ of (\ref{Riccati}) satisfying
(\ref{5.2})--(\ref{5.6}) is called a {\it regular} solution of
(\ref{Riccati}). The following result shows that the regular
solution of (\ref{Riccati}) is unique.

\ms

\bf Corollary 5.3. \sl Let {\rm(SG1)--(SG2)} hold. Then the Riccati
equation (\ref{Riccati}) admits at most one solution $P(\cd)\in
C([t,T];\dbS^n)$ such that (\ref{5.2})--(\ref{5.6}) hold.

\ms

\it Proof. \rm Consider Problem $\hb{(SG)}^0$. Then the adapted solution
$(\eta(\cd),\z(\cd))$ of (\ref{5.8}) is $(\eta(\cd),\z(\cd))\equiv(0,0)$.
Suppose that $P(\cd)$ and $\bar{P}(\cd)$ are two solutions of Riccati equation
(\ref{Riccati}) satisfying (\ref{5.2})--(\ref{5.6}). By Theorem 5.2, we have
$$\lan P(t)x,x\ran=2V(t,x)=\lan \bar{P}(t)x,x\ran,\q\forall x\in \dbR,$$
which implies $P(t)=\bar{P}(t)$. By considering  Problem
$\hb{(SG)}^0$ on $[s,T],t<s<T$, we obtain
$$P(s)=\bar{P}(s),\qq\forall s\in [t,T].$$
This proves our claim. \endpf

\section{Some Examples}

We give some examples in this section. Example 6.1 shows that the
solvability of the Riccati differential equation is not sufficient
enough for the existence of a closed-loop saddle point. So the
$L^2$-regularity is necessary; From Example 6.2 we can see that the
solution of the Riccati equation may be non-unique; To illustrate
that the existence of a closed-loop saddle point does not imply the
existence of an open-loop saddle point, we present Example 6.3.

\ms

\bf Example 6.1. \rm Consider the following optimal control (one-player game) problem:
$$\left\{\2n\ba{ll}
\ns\ds dX(s)=u(s)ds+u(s)dW(s),\q s\in[0,1],\\
\ns\ds X(0)=x,\ea\right.$$
with the cost functional
$$J(x;u(\cd))={1\over2}\dbE X^2(1)+{1\over2}\dbE\int^1_0\big({1\over2}s^3-s^2\big)u^2(s)ds.$$
In this example,
$$\left\{\2n\ba{ll}
\ns\ds A=0,\q B=1,\q b=0,\q C=0,\q D=1,\q \si=0,\\
\ns\ds G=1,\q g=0,\q Q=0,\q\ S=0,\q R(s)={1\over2}s^3-s^2,\q q=0,\q\
\rho=0.\ea\right.$$
The corresponding Riccati equation reads
\bel{6.1}\left\{\2n\ba{ll}
\ns\ds\dot P(s)={2P^2(s)\over {s^3-2s^2+2P(s)}},\q\ae s\in[0,1],\\
\ns\ds P(1)=1.\ea\right.\ee
It is easy to see that $P(s)=s^2$ is the unique solution of
(\ref{6.1}), and
$$\ba{ll}
\ns\ds B(s)^TP(s)+D(s)^TP(s)C(s)+S(s)=s^2,\\
\ns\ds R(s)+D(s)^TP(s)D(s)={1\over2}s^3\ges0,\ea\qq s\in[0,1].$$
Thus, (\ref{5.2}) holds. Now, if the problem has a closed-loop
optimal control, then we should have
$$\Th^*(s)=-\big[R(s)+P(s)\big]^{-1}P(s)=-{2\over s},\qq s\in(0,1],$$
which is not in $L^2(0,1;\dbR)$. This means that the problem does not
have a closed-loop optimal control.

\ms

\bf Example 6.2. \rm Consider the following one-dimensional
controlled system
$$\left\{\2n\ba{ll}
\ns\ds dX(s)=\big[A(s)X(s)+B(s)u(s)\big]ds
+u(s)dW(s),\q s\in[0,1],\\
\ns\ds X(0)=x,\ea\right.$$
with cost functional
$$J(x;u(\cd))={1\over2}\dbE\Big[-X(1)^2+\int_0^1\(Q(s)X(s)^2+R(s)u(s)^2\)ds\Big],$$
where
$$A={1\over2}\Big[{(R-1)^2\over R^2}-1\Big],\q B={R-1\over R},
\q Q=-{1\over R},\q R(s)=\(s-{3\over2}\)^2+{3\over4}>0.$$
The corresponding Riccati equation reads:
\bel{6.3}\left\{\2n\ba{ll}
\ns\ds\dot P+2AP+Q-{B^2P^2\over R+P}=0,\q \ae s\in[0,1],\\
\ns\ds P(1)=-1.\ea\right.\ee
Note that
$$B^2-2A=1,\qq Q+2AR=-2,\qq QR=-1.$$
Then,
$$\ba{ll}
\ns\ds{B^2P^2\over R+P}-2AP-Q={(B^2-2A)P^2-(Q+2AR)P-QR\over
R+P}={P^2+2P+1\over R+P},\ea$$
and (\ref{6.3}) becomes
\bel{6.4}\left\{\2n\ba{ll}
\ns\ds\dot P(s)={P(s)^2+2P(s)+1\over R(s)+P(s)},\q\ae s\in[0,1],\\
\ns\ds P(1)=-1,\ea\right.\ee
which has two solutions:
$$P_1(s)=-1,\qq s\in[0,1],$$
and
$$P_2(s)=s-2,\qq s\in[0,1].$$
We have
$$R(s)+P_1(s)=s^2-3s+2=(s-1)(s-2)\ges0,\qq s\in[0,1],$$
and
$$R(s)+P_2(s)=s^2-2s+1=(s-1)^2\ges0,\qq s\in[0,1].$$
Now, we have
$$\ba{ll}
\ns\ds 2J(x;u(\cd))=\dbE\Big[-X(1)^2+\int_0^1
\(Q(s)X(s)^2+R(s)u(s)^2\)ds\Big]\\
\ns\ds\qq\qq=P(0)x^2+\dbE\int_0^1\Big[\(\dot
P(s)+2A(s)P(s)+Q(s)\)X(s)^2\\
\ns\ds\qq\qq\qq\qq\qq+2P(s)B(s)X(s)u(s)+[R(s)+P(s)]u(s)^2\Big]ds\\
\ns\ds\qq\qq=P(0)x^2+\dbE\int_0^1[R(s)+P(s)]\Big|u(s)+{B(s)P(s)\over
R(s)+P(s)}X(s)\Big|^2ds\\
\ns\ds\qq\qq=P(0)x^2+\dbE\int_0^1[R(s)+P(s)]\Big|u(s)+{[R(s)-1]P(s)\over
R(s)[R(s)+P(s)]}X(s)\Big|^2ds.\ea$$
Note that
$${(R-1)P_1\over R(R+P_1)}={(R-1)(-1)\over R(R-1)}=-{1\over R},$$
and
$${(R-1)P_2\over R(R+P_2)}={(s^2-3s+2)(s-2)\over(s^2-3s+3)(s^2-2s+1)}
={(s-2)^2\over(s^2-3s+3)(s-1)}.$$
Thus,
$$\ba{ll}
\ns\ds
2J(x;u(\cd))=-x^2+\dbE\int_0^1(s-1)(s-2)\Big|u(s)-{X(s)\over(s-1)(s-2)+1}\Big|^2ds\\
\ns\ds\qq\qq\q\ges-x^2=2J(x;u^*(\cd)),\ea$$
with
$$u^*(s)={X(s)\over(s-1)(s-2)+1}\equiv{X(s)\over R(s)}=-{B(s)P_1(s)\over R(s)+P_1(s)}X(s),\q s\in[0,1],$$
which is an optimal control. The closed-loop system reads
$$\left\{\2n\ba{ll}
\ns\ds dX=\Big[{1\over2}\({(R-1)^2\over R^2}-1\)
+\({R-1\over R^2}\)\Big]Xds+{1\over R}XdW(s),\q s\in [0,1],\\
\ns\ds X(0)=x,\ea\right.$$
which is well-posed. Thus, optimal control exists, but Riccati
equation (\ref{6.4}) has more than one solution.

\ms

On the other hand, by taking $P(s)=P_2(s)=s-2$, we have
$$J(x;u(\cd))=-x^2+{1\over2}\dbE\int_0^1(s-1)^2
\Big|u(s)+{(s-2)^2\over(s^2-3s+3)(s-1)}X(s)\Big|^2ds.$$
If
$$\bar u(s)=-{(s-2)^2\over(s^2-3s+3)(s-1)}X(s),$$
is an optimal control, the closed-loop system reads
$$\left\{\2n\ba{ll}
\ns\ds dX=\Big[{1\over2}\({(R-1)^2\over R^2}-1\)
-\({R-1\over R}\){(s-2)^2\over(s^2-3s+3)(s-1)}\Big]Xds\\
\ns\ds\qq\qq\qq\qq\q-\Big[{(s-2)^2\over(s^2-3s+3)(s-1)}\Big]XdW(s),\q s\in [0,1],\\
\ns\ds X(0)=x,\ea\right.$$
which is {\it not} well-posed, since
$$\bar\Th(s)\equiv-{(s-2)^2\over(s^2-3s+3)(s-1)}\not\in L^2(0,1;\dbR).$$
Thus, $\bar u(\cd)$ is not an optimal control.

\ms

\bf Example 6.3. \rm Consider the following state equation
\bel{6.5}\left\{\2n\ba{ll}
\ns\ds dX(s)=\big[u_1(s)-u_2(s)\big]ds+\big[u_1(s)-u_2(s)\big]dW(s),\q s\in[0,1],\\
\ns\ds X(0)=x,\ea\right.\ee
and performance functional
$$J(x;u_1(\cd),u_2(\cd))={1\over2}\dbE\Big[|X(1)|^2+\int_0^1\big(|u_1(s)|^2-|u_2(s)|^2\big)ds\Big].$$
The corresponding Riccati equation reads
$$\left\{\2n\ba{ll}
\ns\ds\dot P=P(1,-1)\begin{pmatrix}1+P&-P\\
-P & -1+P\end{pmatrix}^{-1}\begin{pmatrix}1\\-1\end{pmatrix}P=0,\q 0\les s\les1,\\
\ns\ds P(1)=1.\ea\right.$$
Clearly, $P(s)\equiv 1$. By Theorem 5.2, $(-1,-1)$ is the unique closed-loop saddle point.

\ms

On the other hand, for any $u_1(\cd)\in L_\dbF^2(0,1; \dbR)$,
taking $u_2(\cd)=u_1(\cd)-\l,\ \l\in \dbR$,
the corresponding solution of (\ref{6.5}) is
$$X(s)=x+\lambda s+\lambda W(s),$$
and
$$\dbE|X(1)|^2=x^2+2\l^2+2\l x,$$
$$\dbE\int_0^1\big(|u_1(s)|^2-|u_2(s)|^2\big)ds=-\l^2+2\l\dbE\int_0^1u_1(s)ds.$$
Letting $\l\rightarrow\infty$, we see that
$$\sup_{u_2(\cd)\in L_\dbF^2(0,1;\dbR)} J(x;u_1(\cd),u_2(\cd))
=\infty,\ \forall\ u_1(\cd)\in L_\dbF^2(0,1;\dbR),$$
and hence
$$V^+(0,x)=\inf_{u_1(\cd)\in L_\dbF^2(0,1;\dbR)}\sup_{u_2(\cd)\in L_\dbF^2(0,1;\dbR)}
J(x;u_1(\cd),u_2(\cd))=\infty.$$
So the open-loop saddle point does not exist.

\ms

\section{Concluding Remarks}

In this paper, we present a characterization of the existence of
closed-loop saddle point of linear-quadratic two-person zero-sum
stochastic differential game in terms of the existence of a solution
to differential Riccati equation, with certain regularity. There are
at least two major problems left untouched: (i) The solvability of
the Riccati equation. We mention here that some relevant results can
be found in \cite{Ait Rami-Moore-Zhou 2001} and \cite{McAsey-Mou
2006}. But more complete results are desirable. (ii) The random
coefficients case. This will lead to more involved issues, for
example, the corresponding Riccati equation should be a BSDE, as
indicated in \cite{Chen-Yong 2000,Chen-Yong 2001} for LQ stochastic
optimal control problems with random coefficients. We hope to report
some results relevant to the above-mentioned problems in our future
publications.


\end{document}